\pdfoutput=1
\documentclass[mnsc]{informs_modified}

\SingleSpacedXI

\usepackage[textsize=tiny]{todonotes}

\newcommand{\cp}[2][color= blue!50!red!30!white,linecolor=red,tickmarkheight=0.1cm,noinline]{\todo[#1]{{\bf Chara:} #2}}

\usepackage[ruled,noend]{algorithm2e}

\usepackage[colorlinks=true,citecolor=blue]{hyperref}
\usepackage{cleveref}
\usepackage{threeparttable} % Ensure this package is included
\usepackage{array} % For new column types

\usepackage{natbib}
 \bibpunct[, ]{(}{)}{,}{a}{}{,}%
 \def\bibfont{\small}%
\usepackage{color}
\usepackage{setspace}
\usepackage{bm}
\usepackage{tabu}
\usepackage{array}
\usepackage{multirow}
\usepackage{comment}
\usepackage[toc,page]{appendix}
\usepackage{appendix}
\usepackage{setspace}
\usepackage{hyperref}

\TheoremsNumberedThrough     % Preferred (Theorem 1,
\ECRepeatTheorems

\EquationsNumberedThrough    % Default: (1), (2), ...

\usepackage{tabularx}
\usepackage{hyperref}
\usepackage{url}
\usepackage{booktabs}
\def\WM{W^{\text{n}}}

\newcommand{\xhdr}[1]{\vspace{1mm} \noindent{\bf #1}}

\newcommand{\vstar}{v^{\star}}
\newcommand{\Dstar}{v^{(\textsc{h})}}

\newcommand{\Dtwo}{\textsf{DynUp-2}\xspace}
\newcommand{\Dn}{\textsf{DynUp-n}\xspace}

\usepackage{enumitem}

\begin{document}

\ARTICLEAUTHORS{
\AUTHOR{Patrick Jaillet}\AFF{Department of Electrical Engineering and Computer Science, Massachusetts Institute of Technology,
\EMAIL{jaillet@mit.edu}}
\AUTHOR{Chara Podimata}\AFF{Sloan School of Management, Massachusetts Institute of Technology,
\EMAIL{podimata@mit.edu}}
\AUTHOR{Andrew Vakhutinsky}\AFF{Oracle Lab,
\EMAIL{andrew.vakhutinsky@oracle.com}}
\AUTHOR{Zijie Zhou}\AFF{Operations Research Center,
Massachusetts Institute of Technology, Cambridge, MA, 02139,
\EMAIL{zhou98@mit.edu}}
}

\RUNAUTHOR{}
\RUNTITLE{Online Upgrading Mechanism}

    %\TITLE{Online Upgrading Mechanisms}
    \TITLE{When Should you Offer an Upgrade: Online Upgrading Mechanisms for Resource Allocation}

\ABSTRACT{In this work, we study an upgrading scheme for online resource allocation problems. We work in a sequential setting, where at each round a request for a resource arrives and the decision-maker has to decide whether to accept it (and thus, offer the resource) or reject it. The resources are ordered in terms of their value. If the decision-maker decides to accept the request, they can offer an \emph{upgrade-for-a-fee} to the next more valuable resource. This fee is dynamically decided based on the currently available resources. After the upgrade-for-a-fee option is presented to the requester, they can either accept it, get upgraded, and pay the additional fee, or reject it and maintain their originally allocated resource. 

We take the perspective of the decision-maker and wish to design upgrading mechanisms in a way that simultaneously maximizes revenue and minimizes underutilization of resources. Both of these desiderata are encapsulated in a notion of \emph{regret} that we define, and according to which we measure our algorithms' performance. We present a fast algorithm that achieves $O(\log T)$ regret. Finally, we implemented our algorithm utilizing data akin to those observed in the hospitality industry and estimated our upgrading mechanism would increase the annual revenue by over $17\%$.}

\maketitle

\section{Introduction}

In online resource allocation, we model the interplay between a \emph{supplier} that allocates resources and \emph{requesters} that arrive sequentially and request resources. In general, the resources can be of variable value both to the supplier and the requester, i.e., some are more premium than others. The supplier is tasked with making real-time, irrevocable decisions concerning the method and pricing of the resource to be allocated, despite lacking complete foresight into future demand. Once allocated, the requester provides some payment to the supplier in return. Online resource allocation models have been studied for many different settings such as airline seat allocation \citep{bertsimas2005simulation,zhang2009pricing} and hotel room allocation \citep{ivanov2012hotel,ivanov2014hotel}; their impact has been tremendous both in the academic and the industry world.    
%The practice of online resource allocation can be segmented into two primary categories: profit-based and non-profit-based. In the case of profit-based resource allocation, businesses such as airlines \citep{bertsimas2005simulation,zhang2009pricing} and hotels \citep{ivanov2012hotel,ivanov2014hotel} allocate products, namely seats and rooms, to customers arriving sequentially. These customers, upon receiving the resources, are then obligated to compensate by paying the established price. Conversely, in non-profit scenarios, organizations may provide essential services, like shelters or food, to the sequentially arriving needy, such as the homeless. As a form of reciprocation, these beneficiaries often engage in volunteer activities benefiting the wider community.

An important challenge in online resource allocation occurs when a supplier does not have information about the future realized demand from the requesters and the availability of the resources remains constant. For example, the number of seats on an airplane or rooms in a hotel is fixed, while demand is subject to variability. Demand can surge during certain seasons and shrink during others, but the supplier lacks the flexibility to adjust its resource availability accordingly. As a consequence, there are cases when premium resources can be underutilized. For example, according to \citet{CNN}, approximately $80\%$ of hotel rooms in the United States remain vacant. Addressing this inefficiency calls for the development of advanced resource allocation algorithms, ensuring that even during periods of low demand for premium resources, they are optimally utilized.

One popular idea for ensuring that premium resources do not go underutilized is to employ some \emph{upgrading} scheme. Existing models in the literature \citep{shumsky2009dynamic,yu2015dynamic,cui2023dynamic} primarily assume that suppliers can only offer upgrades to requesters at no extra charge. This assumption (although helpful in simplifying the mathematical modeling) diverges from what happens in reality, where upgrades come at an \emph{extra cost}. For instance, when booking a hotel room through an online platform, it is common for a requester to encounter offers for upgrading to a premium room at an additional fee, typically highlighted with prompts such as "Would you like to upgrade to a premium room for an additional \$50?" This reflects the reality that upgrades are not always complimentary, and determining the optimal pricing for such upgrades is crucial for effective revenue management. %\cpcomment{TODO: Jerry, add a paragraph here that explains what other people in the literature did in terms of the upgrading models they consider}{\color{blue}Jerry:Added}

\subsection{Our Contributions and Techniques}

\xhdr{Model.} We subscribe to the agenda of using upgrades as a way of using premium resources properly and introduce an online upgrading mechanism \emph{for a fee}. At a high level, our model is as follows. Resources of interest (e.g., airline/railway seats, hotel rooms, ride-hailing services, gaming equipment etc.) exhibit an ascending order in both quality and cost. The online upgrading mechanism operates by sequentially processing resource requests. Decisions to accept or reject these requests, along with potential upgrade plans, are made based on factors like arrival rates, remaining availabilities, and associated costs. An upgrade plan entails offering the customer an option to move to the next superior resource category at an additional price. Customers can either accept the upgrade, incurring the extra cost for better resources, or continue using their initially chosen resource at the original price. Our formal model is given in Section~\ref{sec:model}.

\xhdr{Upgrading Mechanism for Online Resource Allocation.} Our main contribution is the development of a mechanism that we call \textsf{DynUp} (stands for ``Dynamic Upgrading'') for identifying \emph{when} and \emph{for how much} a supplier should offer an upgrade. For ease of exposition, we first explain the mechanism \textsf{DynUp-2} and its analysis for the simpler case of $2$ resource types in Section~\ref{sec:resource2}. We use \textsf{DynUp-2} as a building block for our general mechanism \textsf{DynUp-n} (presented in Section~\ref{sec:multiple}) which works for $n$ resource types. The performance of our algorithms is measured in terms of \emph{regret} to an in-hindsight optimum benchmark. Roughly, our algorithms achieve the following guarantees. 

\begin{theorem}[Informal]
    For online resource allocation with $2$ resource types, \textsf{DynUp-2} incurs regret $O(\log T)$. For online resource allocation with $n$ resource types, \Dn incurs regret $O(n \log T)$.
\end{theorem}

Before we delve into the high level description of our mechanism for the $2$-type resource case, we first need to explain the benchmark we are comparing against. At first sight, it may seem tricky to establish a reliable benchmark for the online upgrading problem since at every period we not only have to decide whether to accept or reject a request but also decide on an upgrade price if we do decide to accept. Importantly, this decision at a period $t$ crucially affects the space of allowable decisions for future periods; this is because if we do offer an upgrade price that is ultimately accepted, then we are taking one unit of the premium resource away from future rounds. To navigate this challenge, we propose a ``hybrid programming'' approach which corresponds to an upper bound of the hindsight optimal expected revenue. The ``hybrid programming'' interpolates between a \emph{deterministic} and a \emph{hindsight optimum} counterpart, which draw intuition from approaches traditionally found in price- and quantity- based revenue management respectively. As it will become clear moving forward, we use the ``hybrid programming'' in the development of algorithms \textsf{DynUp-2} and \textsf{DynUp-n}. Detailed explanations of the hybrid programming and the benchmark for $2$ and $n$ resource types can be found in Sections~\ref{subsec:upper} and~\ref{subsec:mupper} respectively.

At a high level, Algorithm \textsf{DynUp-2} determines the upgrading price at each time period $t$ by solving the aforementioned hybrid programming with the remaining availability of each resource and remaining expected demand for each type of resource at $t$. Although the hybrid programming is a non-convex optimization problem, we show that we can obtain a closed-form solution to it. To understand the $O(\log T)$ regret of \Dtwo, we need to study how the upgrade pricing decisions of \Dtwo evolve over time. We distinguish between two cases. First, when the initial resource availability exceeds the expected demand, the stochastic process of \Dtwo's pricing decisions forms a martingale with an expected value change of $0$. Therefore, by comparing this martingale process at each period and the optimal decision variable of ``hybrid programming'', from a known result of \citet{jasin2012re}, we can obtain that the regret can be upper bounded by $O(\log T)$. Second, when the resource availability falls short of expected demand, \Dtwo's dynamic pricing decision is not a martingale. To address this bottleneck we construct a complicated stochastic process which can be shown to be a martingale, closely mirroring the pattern of the actual dynamic pricing decision made by \Dtwo so that the distance of this martingale and the dynamic pricing decision made by \Dtwo at each period can be upper bounded by $O(1/T^2)$.% With this construction, we can also conclude a $O(\log T)$ regret for Algorithm \Dtwo by comparing this constructed martingale process and the optimal decision variable of ``hybrid programming''.

We then extend our findings to the scenario involving $n$ types of resources. Algorithm \Dn is implemented by incorporating Algorithm \Dtwo between consecutive resource types, i.e., between type $i$ and type $i+1$, for all $i \in [n-1]$. In Section \ref{sec:multiple}, we show that the online upgrading problem with $n$ types of resources can be fully decomposed to $n-1$ ``parallel'' problems with $2$ types of resources, leading to a regret bound of $O(n \log T)$. Practically, the number of resource types does not increase with time $T$, making $n$ a fixed parameter. So, under the assumption that $n$ is constant, Algorithm \Dn achieves a regret bound of $O(\log T)$, demonstrating its efficiency and scalability for managing resources across various types.

\xhdr{Empirical Study in Hotel Management.} Our empirical investigation involved a real dataset from the hospitality industry. This dataset encompassed 13,155 requests for rooms spanning one year. Applying Algorithm \Dn to this real-world data led to a significant increase in annual revenue by $17\%$. Notably, in the high-demand month of November, the implementation of Algorithm \Dn resulted in a revenue increase exceeding $26\%$.

\iffalse
\begin{itemize}
    \item Suggest a new mechanism and model.
    \item $O(\log T)$ fast algorithm.
    \item Empirical impact on profit-based resource allocation (Oracle)
    \item Empirical impact on non-profit based resource allocation (need to find a non-profit organization)
\end{itemize}

\fi

\subsection{Related Work} 

\xhdr{Price-based Revenue Management.} \citet{gallego1994optimal, gallego1997multiproduct} introduced the price-based revenue management model. This model processes requests sequentially, with the central decision revolving around the pricing strategy. Typically, higher pricing leads to reduced demand. The authors provided an optimal static pricing technique that achieves an $O(\sqrt{T})$ regret. Subsequently, \citet{jasin2012re, jasin2014reoptimization} developed a dynamic pricing strategy that entails solving a deterministic program for each incoming request, demonstrating that this approach results in an $O(\log T)$ regret. Further advancing this field, \cite{wang2022constant} employed an innovative technique to establish that the dynamic pricing method attains an $O(1)$ regret. Additionally, there are variations to these models, such as \cite{aydin2008pricing}, which explored dynamic pricing in the context of product bundling. 

While our online upgrading process shares certain similarities with price-based revenue management, there are two key distinctions to consider. First, from a modeling standpoint, our approach requires a decision not just on pricing but also on whether to accept each request. This dual consideration positions our model as a hybrid of quantity-based \citep{reiman2008asymptotically,bumpensanti2020re,sun2020near} and price-based models. Second, from technical standpoint, all literature on price-based revenue management assumes a concave revenue function, leading to an online convex optimization problem. In contrast, our model, incorporating both quantity and price variables, is characterized by a non-concave revenue function. %Third, from a societal perspective, dynamic pricing often results in unequal treatment and can even be deemed illegal in some jurisdictions. Our upgrading mechanism, however, presents the same cost to every consumer, thereby upholding a standard of fairness.

\xhdr{Upgrading Models.} The domain of upgrading models has seen various contributions. \citet{shumsky2009dynamic} proposed a dynamic model where upgrades are complimentary, allowing each request to be upgraded only to the immediately superior resource category. \citet{mccaffrey2016optimal} employed dynamic programming to effectively resolve the issue of free upgrading in a two-type resource scenario. Extending this concept, \citet{yu2015dynamic} introduced a multi-step free upgrading framework, permitting requests to be upgraded to any higher-tier resource without cost. \citet{cui2023dynamic} further expanded this model by transitioning from reactive to proactive upgrading. These studies, however, diverge significantly from our work as they are predicated on the assumption of costless upgrading, aligning them more with quantity-based revenue management. In contrast, \citet{gallego2009upgrades} introduced the concept of ``upsell'' as a non-complimentary upgrade in their research as an extension. They focused on static upselling strategies but did not provide any theoretical performance guarantees for this method.

\iffalse
Considerations that have been discussed in the literature: 
\begin{enumerate}
    \item Quantity-based (Involuntary Upgrade) 
    \begin{itemize}
    \item ``Limited Cascading Model'' optimal upgrading, where you can only upgrade to the next high fare class. See e.g., \cite{gallego2009upgrades}, \cite{shumsky2009dynamic}
    \item Optimal upgrading strategy (DP) in 2 dimensions (only business and economy class). See, e.g. \cite{mccaffrey2016optimal}.
    \item ``General'' upgrading, where you don't need to upgrade to the next best resource, but you can upgrade to any better resource (with different costs, of course). See e.g., \cite{yu2015dynamic}
    \item Upgrading on network RM. Customers can request a set of products together. See e.g., \cite{cui2023dynamic}
\end{itemize}
    %\item Pricing
    %\begin{itemize}
    %    \item Static Pricing. See e.g. \cite{gallego1994optimal}, \cite{gallego1997multiproduct}.
    %    \item Dynamic Pricing. See e.g. \cite{jasin2012re}, \cite{jasin2014reoptimization}, \cite{wang2022constant}
    %    \item Dynamic Pricing with Bundling. See e.g., \cite{aydin2008pricing}
    %\end{itemize}
\end{enumerate}
\fi

\section{Model}\label{sec:model}
We consider a setting where a supplier wants to allocate $n$ distinct resources to interested individuals (aka \emph{requesters}). Each resource has a different type $i \in [n]$.
For each resource type $i \in [n]$, there is an associated availability $c_i$ and a cost $r_i$. We assume that resources are ordered in terms of their cost to the requesters and that higher costs are associated with higher value for all the requesters, i.e., $r_1<r_2<\ldots<r_n$ and the $n^{\text{th}}$ resource type is the most premium. On the request side, requests arise sequentially across a discrete time horizon spanning $T$ periods. In any given period, at most one requester may show up, and request precisely one unit of some resource. The request arrival rate for resource type $i$ is known and denoted by $\lambda_i$, and $\lambda_0$ represents the rate of no request. We normalize the arrival rates such that $\sum_{i=0}^{n} \lambda_i = 1$. Consequently, the probability of a request for resource type $i$ within a specific period is $\lambda_i$, and the probability for no request is $\lambda_0$.

When a request for a resource of type $i$ emerges, the supplier is faced with the decision to either accept or reject this request. If the request is rejected, there are no resources allocated. On the other hand, if the request is accepted, the supplier offers an upgrade plan\footnote{If $i=n$ or the type $i+1$ resource has availability less than $1$, no upgrade plan will be offered.}: would the requester be willing to pay an additional cost, $u_i \in [0,r_{i+1}-r_i]$, in exchange for the superior type $i+1$ resource\footnote{If we accept the type $i$ request, but type $i$ resource has availability less than $1$, we can only upgrade this request to the type $i+1$ resource with $u_i=0$. }? If this plan is accepted by the requester, the supplier assigns one unit of the type $i+1$ resource, leading the requester to pay a combined cost of $r_i+u_i$. Conversely, a rejection results in the allocation of the type $i$ resource, with a corresponding cost of $r_i$. We posit that any request for a resource of type $i$ has a $f_i(u_i)$ likelihood of accepting the upgrade for an additional cost of $u_i$. The overall objective for the supplier is to \emph{maximize the revenue} obtained given the resources and their availabilities, and the designated upgrading offer mechanism. We use \emph{regret} to measure the performance of the online upgrading mechanism, where the regret is defined as the hindsight optimal revenue $\mathbb{E} [ W^{\text{OPT}}]$ minus the expected cumulative revenue generated by the mechanism: 
\[
R(T) = \mathbb{E} [ W^{\text{OPT}} - W^\pi]
\]
where $W^\pi$ denotes the revenue generated by a mechanism $\pi$.

\subsection{Model Assumptions and Notations}
We assume that the likelihood functions \( f_i(\cdot) \) are known both by the supplier and the requesters. These functions are bijective, continuous, and monotone-decreasing. Further, let \( f_i(0) = 1 \), capturing the fact that no individual would decline a free
%\cp{The word "involuntary" may suggest that the individual didn't want the upgrade. {\color{blue}Jerry:Maybe I can change the word to free upgrade. In my internship, I think they called free upgrade by involuntary upgrade}} \av{"involuntary" was related to the hotel operator (supplier) in the situation when they did not have any choice but to upgrade a guest for free.  It's a kind of industry slang and probably should be avoided here} 
upgrade. The inverse function of \( f_i(\cdot) \) is denoted by \( p_i(\cdot) \) (e.g., if \( v_i = f_i(u_i) \), then \( u_i = p_i(v_i) \)). In other words, function $p_i(v)$ returns the upgrade price $u$ such that the probability that $u$ is accepted for the upgrading from resource $i$ to $i+1$ is $f_i(u)$. Furthermore, the functions \( p_i(\cdot) \) are also bijective, continuous, and monotone decreasing. Since $f_i(0) =1$, then \( p_i(1) = 0 \). The expected revenue arising from the upgrade are captured by \( R_i(v_i) = v_i p_i(v_i) = u_i f_i(u_i) \). %Aligning with prevailing conventions in dynamic pricing literature \citep{gallego1994optimal,gallego1997multiproduct,jasin2012re,jasin2014reoptimization}, 
We assume that the functions \( R_i(\cdot) \) are bounded and quasi-concave. The assumption of quasi-concavity aligns with common functional forms of $p_i(\cdot)$, such as linearly and exponentially decreasing functions, which inherently lead to quasi-concave $R_i(\cdot)$. We also assume that $R_i(\cdot)$ has a finite maximizer, denoted by \( \vstar_i \), which lies in the interval \([ f_i(r_{i+1} - r_i), 1 ]\).

%Question: what is the academic words for generalized resource allocation? Specifically,
%\begin{itemize}
%    \item Total revenue $\to$ (revenue? (value of resource allocation))
%    \item Price $\to$  (utility or value)
%    \item Payment of Each request $\to$ (cost)
%\end{itemize}

% “total revenue” instead of “total revenue”, “social cost” instead of “price”, and “social benefit” instead of “pay the amount”.

\section{Algorithm for Two Types of Resources}\label{sec:resource2}
In this section, we focus on the setting with two distinct resource types where only one of them can be upgraded to the other. We assume that the availability of the premium resource (type 2) exceeds its expected demand, namely \( c_2 > \lambda_2 T \). Otherwise, the problem becomes trivial as there would be no need for upgrading. We first describe the optimal stochastic and deterministic formulations of the online upgrading mechanism  to introduce the upper bound performance for any upgrading algorithm (Section \ref{subsec:upper}). Second, we present our algorithm, and its regret analysis in Section \ref{subsec:alg2}.

\subsection{Benchmark Revenue Upper Bound} \label{subsec:upper}

In each period \( t \in [T] \), let $\Lambda_i^{(t)}$ be the binary variable that represents whether the type $i$ request arrives in period $t$, and \( y_i^{(t)} \) be the binary variable that indicates whether the supplier accepts the type \( i \) request. Let \( D_i^{(t)} \) represent the quantity of type \( i \) resources consumed in period \( t \). For example, at period $t$, if a type $1$ request arrives and accepts the upgrading plan, then no type $1$ resource is used and one unit of type $2$ resource is consumed, and hence $D_1^{(t)}=0$ and $D_2^{(t)}=1$. Note that \( D_i^{(t)} \) is a random variable, influenced by the upgrading cost \( u_i \). Hence, the optimal cumulative revenue can be expressed as the following stochastic programming:
\begin{align} \label{eq:optdef}
    W^{\textsc{opt}} \nonumber& = \max_{\mathbf{y_1},\mathbf{y_2},\mathbf{u_1}} \sum_{t\in [T]} \mathbb{E} \left[\Lambda_1^{(t)}y_1^{(t)}\left(r_1 D_1^{(t)} + (r_1+u_1^{(t)})D_2^{(t)} \right) + \Lambda_2^{(t)} y_2^{(t)}r_2D_2^{(t)}  \right] \\ \nonumber& \text{s.t. } \sum_{t\in [T]} \Lambda_1^{(t)}y_1^{(t)}D_1^{(t)} \leq c_1
    \\& \text{ } \text{ } \text{ } \text{ } \text{ } \sum_{t \in [T]} (\Lambda_1^{(t)}y_1^{(t)}+\Lambda_2^{(t)}y_2^{(t)})D_2^{(t)} \leq c_2,   \tag{SP}
\end{align}
where the expectation is taken on the randomness of the arrival process, and the constraints must hold almost surely. In (SP), we maximize the sum of the expected revenue over all periods $t \in [T]$. In each period $t$, if a type $1$ request shows up ($\Lambda_1^{(t)}=1$) and is accepted ($y_1^{(t)}=1$), the expected revenue generated is $r_1D_1^{(t)}$ (the request rejects the upgrading plan) plus $(r_1+u_1^{(t)})D_2^{(t)}$ (the request accepts the upgrading plan). If a type $2$ request shows up ($\Lambda_2^{(t)}=1$) and is accepted ($y_2^{(t)}=1$), the expected revenue generated is $r_2D_2^{(t)}$. The first constraint guarantees that the sum of the type $1$ resource consumed over all periods does not exceed $c_1$, and the second constraint ensures that the sum of the type $2$ resource consumed (type $1$ request accepts the upgrading plan or type $2$ request is accepted) over all periods does not surpass $c_2$. 

Observe that \ref{eq:optdef} only outlines the decision-making process under the condition that both resource types remain available; it does not account for scenarios where either resource is exhausted. For example, in cases where type $2$ resource runs out, there is no upgrading option. Under these circumstances, the approach mandates rejecting all type $2$ requests, while accepting those from type $1$, provided there is remaining availability in type $1$ resources. To solve the \eqref{eq:optdef} we can use the following dynamic program: \begin{align*}
W_t(c_1^{(t)},c_2^{(t)}) = \max_{y_1^{(t)},y_2^{(t)},u_1^{(t)}}\mathbb{E} & \Big[\Lambda_1^{(t)}y_1^{(t)}\left(r_1 D_1^{(t)} + (r_1+u_1^{(t)})D_2^{(t)} \right) + \Lambda_2^{(t)}y_2^{(t)}r_2D_2^{(t)} \\& + W_{t+1}(c_1^{(t)}-\Lambda_1^{(t)}y_1^{(t)}D_1^{(t)},c_2^{(t)}-(\Lambda_1^{(t)}y_1^{(t)}+\Lambda_2^{(t)}y_2^{(t)})D_2^{(t)})  \Big]. \end{align*}
where $c_1^{(t)},c_2^{(t)}$ are the remaining availability for type $1$, $2$ resource respectively at time $t$. However, the curse of dimensionality presents a significant challenge to computing the exact solution of the dynamic program. To counteract this problem, a prevalent approach in the context of price-based revenue management is to formulate a \emph{deterministic} counterpart to the stochastic program defined in \eqref{eq:optdef}~\citep{gallego1994optimal, gallego1997multiproduct}. This is achieved by substituting the upgrade price random variable $v_1^{(t)}=f_1(u_1^{(t)})$ with a fixed variable $v_1$. However, contrary to price-based revenue management, in our problem, there also exists a \emph{quantity-based} variable $y_i^{(t)}$ capturing whether to accept or reject the arriving type $i$ request and we need to account for it. A well-known formulation that can be treated as an upper bound of the (SP) and contains quantity-based variables is the \emph{hindsight optimum formulation} \citep{reiman2008asymptotically,bumpensanti2020re}, where we define $Y_i$ as the total number of type $i$ requests to accept, and define $\Lambda_i$ as the aggregate number of arrivals for each type $i$. Here, $\Lambda_i$ is a binomial random variable, namely $\Lambda_i \sim \text{Bin}(T,\lambda_i)$, and in the hindsight formulation, one assumes that the value of $\Lambda_i$ is known in advance.

Putting everything together, we propose a \emph{hybrid programming} that integrates both the \emph{hindsight} and \emph{deterministic} formulations, where $\Lambda_i$ is the hindsight variable and $v_1$ is the deterministic variable. This formulation will henceforth be referred to as ``hybrid programming'' (HP).
\begin{align} \label{eq:optD}
    w^{\textsc{h}} \nonumber& = \max_{Y_1,v_1} \Lambda_2r_2 + Y_1\left( (r_1+p(v_1))v_1+r_1(1-v_1) \right) \\ \nonumber& \text{s.t. } Y_1 \leq \Lambda_1
    \\ \nonumber& \text{ } \text{ } \text{ } \text{ } \text{ } Y_1v_1  \leq c_2-\Lambda_2 \\& \text{ } \text{ } \text{ } \text{ } \text{ } Y_1(1-v_1) \leq c_1.
\end{align}
We provide some more details about the optimization problem in~\eqref{eq:optD}. The objective function captures the total revenue. Since we always accept type $2$ requests when there is availability, the revenue generated by type $2$ requests is $\Lambda_2r_2$. For type $1$ requests, we accept $Y_1$ of them and each request will generate an expected revenue of $(r_1+p(v_1))v_1$ (accept the upgrading plan) plus $r_1(1-v_1)$ (reject the upgrading plan). The first constraint ensures that the number of type $1$ requests accepted is no larger than the total number of type $1$ arrivals. The second constraint guarantees that the total number of upgraded type $1$ requests is less than or equal to the availability of type $2$ resource minus the number of type $2$ requests. The third constraint suggests that the total number of type $1$ requests who reject the upgrading plan is less than or equal to the availability of type $1$ resource. As $(r_1+p(v_1))v_1+r_1(1-v_1) $ is always positive, the optimization problem \eqref{eq:optD} is equivalent to the following hybrid of hindsight and deterministic programming:
\begin{equation} \label{eq:optD1}
     w^{\textsc{h}} = \max_{v_1} \Lambda_2r_2 + \min\left\{\Lambda_1, \frac{c_2-\Lambda_2}{v_1}, \frac{c_1}{1-v_1} \right\}\left( R_1(v_1)+r_1 \right) \tag{HP}.
\end{equation}

Let $w^{\textsc{h}}$ be the optimal objective value; note that $w^{\textsc{h}}$ is a random variable that depends on $\Lambda_1$ and $\Lambda_2$. We define $W^{\textsc{h}}=\mathbb{E}[w^{\textsc{h}}]$.

Recall that \eqref{eq:optdef} only captures the decision-process only up to the point where a resource runs out. Consequently, $W^{\textsc{h}}$ should be regarded as an approximate upper bound for the total revenue achievable prior to the depletion of any resources.  In our framework, decision-making ends either at the exhaustion of both resources or upon surpassing the designated time horizon, denoted as $[0, T]$. Accordingly, it becomes necessary to split the analysis based on the order of resource depletion, with comprehensive elaborations provided in the following theorem.

%In contrast to typical findings in revenue management, $W^{\textsc{h}}$ is \textbf{not} an upper bound for \eqref{eq:optdef}. %\cpcomment{i don't understand the point that you are making here. I think it is subtle and I would like for us to explain it more clearly. Specifically, I thought that $w^H$ was formulated as relaxations or upper bounds to (SP). But it seems there is something subtle here that has to do with the order of the resource depletions...?} 
%This deviation arises because, in standard revenue management models, decisions stop upon the depletion of any resource. However, in our framework, decision-making ends either at the exhaustion of both resources or upon surpassing the designated time horizon, denoted as $[0, T]$. Accordingly, it becomes necessary to split the analysis based on the order of resource depletion, with comprehensive elaborations provided in the following theorem. 

\begin{theorem} \label{thm:sameopt}
    Define three sets $\mathcal V_1$, $\mathcal V_2$, $\mathcal V_3$ such that $\mathcal V_1 = \left\{v_1: \Lambda_1 = \min\{\Lambda_1, \frac{c_2-\Lambda_2}{v_1}, \frac{c_1}{1-v_1} \} \right\}$, $\mathcal V_2 = \left\{v_1: \frac{c_2-\Lambda_2}{v_1} = \min\{\Lambda_1, \frac{c_2-\Lambda_2}{v_1}, \frac{c_1}{1-v_1} \} \right\}$, and $\mathcal V_3 = \left\{v_1: \frac{c_1}{1-v_1} = \min\{\Lambda_1, \frac{c_2-\Lambda_2}{v_1}, \frac{c_1}{1-v_1} \} \right\}$. 
        \begin{align*}
        \text{For $v_1 \in \mathcal V_1$:} &\qquad
        &&w^{\textsc{U}_1}= \max_{v_1 \in \mathcal V_1} \Lambda_2r_2 + \Lambda_1 \left( R_1(v_1)+r_1 \right) \\ 
        \text{For $v_1 \in \mathcal V_2$:} &\quad
        &&w^{\textsc{U}_2} = \max_{v_1 \in \mathcal V_2} \Lambda_2r_2 +\frac{c_2-\Lambda_2}{v_1} R_1(v_1) + r_1\min\left\{\Lambda_1, c_1+c_2-\Lambda_2 \right\}\\
        \text{For $v_1 \in \mathcal V_3$:} &\quad &&w^{\textsc{U}_3}= \max_{v_1 \in \mathcal V_3} \Lambda_2r_2 + \frac{c_1}{1-v_1} R_1(v_1) + r_1\min\left\{\Lambda_1, c_1+c_2-\Lambda_2 \right\}.
        \end{align*}
    % \begin{itemize}
    %     \item For $v_1 \in \mathcal V_1$,
    %     \[
    %     w^{\textsc{U}_1} = \max_{v_1 \in \mathcal V_1} \Lambda_2r_2 + \Lambda_1 \left( R_1(v_1)+r_1 \right).
    %     \]
    %     \item For $v_1 \in \mathcal V_2$,
    %     \[
    %     w^{\textsc{U}_2} = \max_{v_1 \in \mathcal V_2} \Lambda_2r_2 +\frac{c_2-\Lambda_2}{v_1} R_1(v_1) + r_1\min\left\{\Lambda_1, c_1+c_2-\Lambda_2 \right\}.
    %     \]
    %     \item For $v_1 \in \mathcal V_3$,
    %     \[
    %     w^{\textsc{U}_3} = \max_{v_1 \in \mathcal V_3} \Lambda_2r_2 + \frac{c_1}{1-v_1} R_1(v_1) + r_1\min\left\{\Lambda_1, c_1+c_2-\Lambda_2 \right\}.
    %     \]
    % \end{itemize}
    Let $w^U$ be such that:
    \begin{equation} \label{eq:upperdp}
    w^{\textsc{U}} = \max\{w^{\textsc{U}_1},w^{\textsc{U}_2},w^{\textsc{U}_3} \}.  \tag{Upper-HP}
    \end{equation}
    Then, given any $\Lambda_1$ and $\Lambda_2$, the optimal solution to \ref{eq:upperdp} is the same as the one to \ref{eq:optD1}. Furthermore, $W^{\textsc{U}} \geq \mathbb{E}[W^{\textsc{opt}}]$, where  $W^{\textsc{U}}=\mathbb{E}[w^{\textsc{U}}]$.
\end{theorem}

The proof of Theorem \ref{thm:sameopt} can be found in Appendix \ref{secappend:upper}. To demonstrate that Equations \eqref{eq:upperdp} and \eqref{eq:optD1} share the same optimal solution, we only need to calculate the optimal solution in three distinct cases: when $v_1$ is inside $\mathcal V_1$, or $\mathcal V_2$, or $\mathcal V_3$. To establish that $w^{\textsc{U}}$ is the maximum potential revenue achievable by any online algorithm, we find that, in each case, upon fixing the hindsight variables, $\Lambda_i$, the optimization problem with respect to $v_1$ is convex. Then, by \cite{gallego1994optimal}, the statement immediately follows. 

Theorem \ref{thm:sameopt} describes hybrid programming formulations \eqref{eq:upperdp} that effectively capture the order of resource depletion. Specifically, $w^{\textsc{U}_1}$ captures the case where no resource will run out because the number of type $1$ arrivals is not enough. $w^{\textsc{U}_2}$ captures the case where the type $2$ resource runs out first, and after that, the supplier cannot give an upgrading option to the remaining type $1$ requests. In this case, the remaining type $1$ requests can only generate a revenue of $r_1$ each. $w^{\textsc{U}_3}$ describes the case where the type $1$ resource first runs out, and after that, the supplier can only give free upgrades to the remaining type $1$ requests. In that case, the remaining type $1$ requests can also only generate a revenue of $r_1$ each. Finally, Theorem \ref{thm:sameopt} states that the optimal solution to \eqref{eq:optD1} coincides with the one to \eqref{eq:upperdp}, namely $\Dstar_1=\text{argmax}_{v_1} w^{\textsc{h}}=\text{argmax}_{v_1} w^{\textsc{U}}$. 

Going back to our benchmark revenue $\mathbb{E}[W^{\text{OPT}}]$, by Theorem \ref{thm:sameopt}, $W^{\textsc{U}}$ is an upper bound for $\mathbb{E}[W^{\textsc{opt}}]$, which implies that  $R(T) \leq \mathbb{E}[W^{\textsc{U}}-W^{\pi}]$, where the expectation is taken on the arrival process. Henceforth, we will use $\mathbb{E}[W^{\textsc{U}}-W^{\pi}]$ as a stronger performance metric. 

% \cpcomment{It is unclear here what you are trying to prove in this subsection. I suggest that you write as precise mathematical statements as possible. At least for EC, we will need the Lemma statements in the main body too. Between each lemma you can have a sentence or two about why it's needed to complete the proof of the result you present in this section.}{\color{blue}Jerry: I have put all lemmas to a single theorem above, and I explain each formulation in the theorem above.}

\subsection{Algorithm Description and Regret Analysis} \label{subsec:alg2}

In this section, we introduce a fast algorithm (Algorithm \Dtwo) for the supplier to make the decision at each time period: either accept or reject the demand, and if they accept, what the upgrade price will be. Formally, we prove the following theorem.
\begin{theorem} \label{thm:2type}
    The regret of Algorithm \Dtwo is bounded by $O(\log T)$.
\end{theorem}

\begin{algorithm}[!tb]
\caption{Optimal Upgrade Price}
\label{alg:price}

\DontPrintSemicolon % To remove semicolons at the end of lines
\KwIn{time period $t$; remaining availability $c_1^{(t)}$, $c_2^{(t)}$; arrival rate $\lambda_1$, $\lambda_2$.}
\KwOut{Optimal upgrade price at time $t$.}
\eIf{$c_1^{(t)} + c_2^{(t)} < (\lambda_1+\lambda_2) (T-t+1)$}{
    Let upgrade plan $v_1^{(t)}=\min \left\{\max \left\{\frac{c_2^{(t)}-\lambda_2(T-t+1)}{c_1^{(t)}+c_2^{(t)}-\lambda_2(T-t+1)} ,0 \right\} ,1 \right\}$\;
}{
    Let upgrade plan $v_1^{(t)}=\max\left\{\min\left\{\vstar_1,\frac{c_2^{(t)}}{\lambda_1(T-t+1)}-\frac{\lambda_2}{\lambda_1}\right\}, 1-\frac{c_1^{(t)}}{\lambda_1(T-t+1)}\right\} $\;
}
\Return{$\min \left\{\max \left\{p_1(v_1^{(t)}),0 \right\} ,1 \right\}$.\;}
\end{algorithm}

\begin{algorithm}[!tb]
\caption{\Dtwo}
\label{alg:2type}

\DontPrintSemicolon % To remove semicolons at the end of lines
\KwIn{time horizon $T$; initial availability $c_1^{(1)}$, $c_2^{(1)}$; arrival rate $\lambda_1$, $\lambda_2$.}
\For{$t \in \{1,2,...,T\}$}{
    Observe demand of type $i$. If $i=2$, we accept the request if and only if $c_2^{(t)} \geq 1$. \;
    \uIf{$i=1$ and $c_1^{(t)} \geq 1$}{
        Input $t$, $c_1^{(t)}$, $c_2^{(t)}$, $\lambda_1$, and $\lambda_2$ to Algorithm \ref{alg:price}, and denote the output as $u_1^{(t)}$.\;
        \uIf{$c_2^{(t)} \geq 1$}{
            Accept the request and give an upgrade plan with additional price $u_1^{(t)}$.\;
        }
        \Else{
            Accept the request without an upgrade plan.\;
        }
    }
    \uElseIf{$i=1$ and $c_1^{(t)} = 0$}{
        Accept the request and give a free upgrade if and only if $\left(c_2^{(t)}-\lambda_2(T-t+1)\right)^{+} > \frac{1}{2}\lambda_1 (T-t+1)$\;
    }
}
\end{algorithm}

At a high level, in Algorithm \Dtwo, we split the time horizon into two segments. The first segment, $[1,\tau]$, starts from the beginning and extends up to the moment $\tau$, defined as the first instance when the type $1$ resource is completely depleted, i.e., $\tau = \inf\{t:c_1^{(t)}=0\}$. The second segment, $[\tau,T]$, covers the remaining portion of the time horizon. In the first segment $[1,\tau]$, requests are accepted on a first-come first-serve basis (i.e., demand is rejected only if there is insufficient availability). The upgrade price is determined by Algorithm \ref{alg:price} in each time step $t$. In the second segment $[\tau,T]$, as the type $1$ resource is depleted, the problem becomes a single-leg revenue management problem: there is one type of resource serving two types of customers with revenue for the supplier $r_1$ and $r_2$ respectively. In this case, the last step of \Dtwo borrows the state-of-the-art algorithm for single-leg revenue management from \citet{vera2021bayesian}. Let $W^{\textsc{U}}(t_1,t_2)$ be the upper bound of the revenue between $[t_1,t_2]$, and $W^{\pi}(t_1,t_2)$ the total revenue generated by Algorithm \Dtwo between $[t_1,t_2]$. From the definition of regret: 
\begin{align}
    R(T) \nonumber&= \mathbb{E}[W^{\textsc{U}}(1,T)-W^{\pi}(1,T)] \\&= \mathbb{E}[W^{\textsc{U}}(1,\tau-1)-W^{\pi}(1,\tau-1)]+ \mathbb{E}[W^{\textsc{U}}(\tau,T)-W^{\pi}(\tau,T)] \label{eq:regr2}%\\&\leq \mathbb{E}\left[\sum_{t=1}^{\tau-1}\left(R(\Dstar_1)-R(v_1^{(t)})\right) \right] + \mathbb{E}[W^{\textsc{U}}(\tau,T)-W^{\pi}(\tau,T)],
\end{align}

Then, we bound the two expectations in Equation \eqref{eq:regr2} respectively. First, we suggest a lemma to upper bound $\mathbb{E}[W^{\textsc{U}}(1,\tau-1)-W^{\pi}(1,\tau-1)]$ as follows:

\begin{lemma}\label{lem:firstupper} The regret of the first time segment is upper bounded as: 
    \begin{center}
    $\mathbb{E}[W^{\textsc{U}}(1,\tau-1)-W^{\pi}(1,\tau-1)] \leq \mathbb{E}\left[\sum_{t \in [\tau-1]}\left(R(\Dstar_1)-R(v_1^{(t)})\right) \right].$
    \end{center}
\end{lemma}

The proof can be found in Appendix \ref{secappend:thm1}. Next, we study the relationship between $\Dstar_1$ and $v_1^{(t)}$. Observe that $v_1^{(t)}$ is a stochastic process, which corresponds to the optimal solution to \eqref{eq:upperdp}, taking into account both the remaining availability of each resource, $c_i^{(t)}$, and the anticipated demand for each category, $\lambda_i(T-t+1)$, $i \in \{1,2\}$. According to Theorem \ref{thm:sameopt},  the optimal solutions for \eqref{eq:upperdp} and \eqref{eq:optD1} are equivalent, and our focus shifts primarily to solving \eqref{eq:optD1}. 
Notice that \eqref{eq:optD1} is a non-convex optimization problem at each period $t$, the following lemma shows the closed-form solutions of \eqref{eq:optD1} when $t=1$.

\begin{lemma} \label{lem:Dstar1}
    If $c_1^{(1)}+c_2^{(1)} > (\lambda_1+\lambda_2)T$, we have that:
    \[
    v_1^{(1)} = \max \left\{ \min\left\{ \vstar_1,\frac{c_2^{(1)}-\lambda_2 T}{\lambda_1 T}   \right\}, \frac{\lambda_1 T-c_1^{(1)}}{\lambda_1 T} \right\}\ \text{and } \, \Dstar_1 = \max \left\{ \min\left\{ \vstar_1,\frac{c_2^{(1)}-\Lambda_2}{\Lambda_1}   \right\}, \frac{\Lambda_1 -c_1^{(1)}}{\Lambda_1 } \right\}.
    \]
    % \[
    % \Dstar_1 = \max \left\{ \min\left\{ \vstar_1,\frac{c_2^{(1)}-\Lambda_2}{\Lambda_1}   \right\}, \frac{\Lambda_1 -c_1^{(1)}}{\Lambda_1 } \right\}.
    % \]
    If $c_1^{(1)}+c_2^{(1)} \leq (\lambda_1+\lambda_2)T$, we have that: 
    \begin{center}
    $v_1^{(1)} = \frac{c_2^{(1)}-\lambda_2T}{c_1^{(1)}+c_2^{(1)}-\lambda_2T} \, \text{and } \, \Dstar_1 = \frac{c_2^{(1)}-\Lambda_2}{c_1^{(1)}+c_2^{(1)}-\Lambda_2}.$
    \end{center}
    % \[
    % \Dstar_1 = \frac{c_2^{(1)}-\Lambda_2}{c_1^{(1)}+c_2^{(1)}-\Lambda_2}.
    % \]
\end{lemma}

\begin{figure}[!tb]
\center
\includegraphics[width=0.75\textwidth]{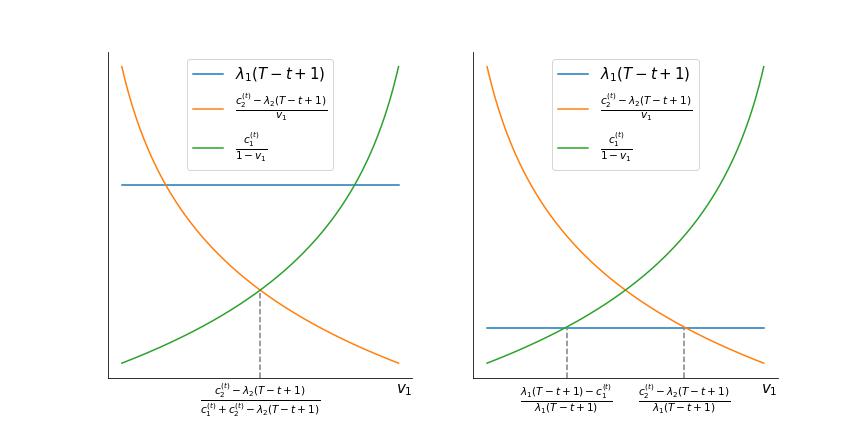}
\caption{Auxiliary image for the proof sketch. {\bf Left}: the case where $c_1^{(t)} + c_2^{(t)} < (\lambda_1+\lambda_2) (T-t+1)$. {\bf Right}: the case where $c_1^{(t)} + c_2^{(t)} \geq (\lambda_1+\lambda_2) (T-t+1)$} 
\label{fig:3function}
\end{figure}

Before providing the full proof, we first give some intuition. In the first case where $c_1^{(t)} + c_2^{(t)} < (\lambda_1+\lambda_2) (T-t+1)$, the aggregate availability is less than the cumulative estimated expected future arrivals; the left diagram in Fig. \ref{fig:3function} shows that $\min\left\{\lambda_1(T-t+1), \frac{c_2^{(t)}-\lambda_2(T-t+1)}{v_1}, \frac{c_1^{(t)}}{1-v_1} \right\}$ cannot equal $\lambda_1(T-t+1)$. If $\frac{c_2^{(t)}-\lambda_2(T-t+1)}{v_1}$ is the minimum, Lemma \ref{lem:monotone} from Appendix \ref{secappend:upper} implies that $W^{\textsc{h}}$ is decreasing with respect to $v_1$. Conversely, if $\frac{c_1^{(t)}}{1-v_1}$ is the minimal value, the same lemma suggests that $W^{\textsc{h}}$ is increasing with respect to $v_1$. Consequently, at time period $t$, the optimal solution to \eqref{eq:optD1} is at the intersection of the functions $\frac{c_2^{(t)}-\lambda_2(T-t+1)}{v_1}$ and $\frac{c_1^{(t)}}{1-v_1}$, yielding $\Dstar_1 = \frac{c_2^{(t)}-\lambda_2(T-t+1)}{c_1^{(t)}+c_2^{(t)}-\lambda_2(T-t+1)}$.

For the second case where $c_1^{(t)} + c_2^{(t)} \geq (\lambda_1+\lambda_2) (T-t+1)$, the right diagram of Fig. \ref{fig:3function} shows the relationship between $\lambda_1(T-t+1)$, $\frac{c_2^{(t)}-\lambda_2(T-t+1)}{v_1}$, and $\frac{c_1^{(t)}}{1-v_1}$. By Lemma \ref{lem:monotone}, the monotonicity of $W^{\textsc{h}}$ in each function shows that the optimal solution $\Dstar_1$ is within 
\begin{center}
    $\mathcal I = \left\{v_1:\lambda_1(T-t+1)=\min\left\{\lambda_1(T-t+1), \frac{c_2^{(t)}-\lambda_2(T-t+1)}{v_1}, \frac{c_1^{(t)}}{1-v_1} \right\} \right\}.$
\end{center}
In this context, \eqref{eq:optD1} becomes a convex optimization problem: $\max_{v_1 \in \mathcal I} R_1(v_1)$. $R_1(v_1)$ has a unique maximizer $\vstar_1$, so we set $v_1^{(t)}$ as the projection of $\vstar_1$ to $\mathcal I$. The formal proof follows.

\proof{Proof of Lemma \ref{lem:Dstar1}}
If $c_1^{(1)}+c_2^{(1)} > (\lambda_1+\lambda_2)T$, Lemma \ref{lem:monotone} implies that $v_1^{(1)} \in \left[\frac{\lambda_1 T-c_1^{(1)}}{\lambda_1 T}, \frac{c_2^{(1)}-\lambda_2 T}{\lambda_1 T}\right]$. Let the unique maximizer of $R_1(v_1)$ be $\vstar_1$. If $\vstar_1 \in \left[\frac{\lambda_1 T-c_1^{(1)}}{\lambda_1 T}, \frac{c_2^{(1)}-\lambda_2 T}{\lambda_1 T}\right]$, then 
\begin{center}
$v_1^{(1)} = \vstar_1 = \text{argmax}_{v_1 \in \left[\frac{\lambda_1 T-c_1^{(1)}}{\lambda_1 T}, \frac{c_2^{(1)}-\lambda_2 T}{\lambda_1 T}\right]} R_1(v_1).$
\end{center}
If $\vstar_1 > \frac{c_2^{(1)}-\lambda_2 T}{\lambda_1 T}$, by Lemma \ref{lem:monotone}, we have $v_1^{(1)} = \frac{c_2^{(1)}-\lambda_2 T}{\lambda_1 T}$. Similarly, if $\vstar_1 < \frac{\lambda_1 T-c_1^{(1)}}{\lambda_1 T}$, by Lemma \ref{lem:monotone}, we have $v_1^{(1)} = \frac{\lambda_1 T-c_1^{(1)}}{\lambda_1 T}$. By the same approach, we can also get the expression of $\Dstar$.

If $c_1^{(1)}+c_2^{(1)} \leq (\lambda_1+\lambda_2)T$, by Theorem \ref{thm:sameopt}, we have in this case, $
W^{\textsc{U}} = \max\{W^{\textsc{U}_2},W^{\textsc{U}_3} \}$. The max can be attained if and only if $W^{\textsc{U}_2}=W^{\textsc{U}_3} $, which implies that 
\[
\frac{c_2-\lambda_2T}{v_1^{(1)}}=\frac{c_1}{1-v_1^{(1)}}.
\]
Because the left hand side is monotone decreasing in $v_1^{(1)}$ and right hand side is monotone increasing in $v_1^{(1)}$, there is only one solution, which is 
\[
v_1^{(1)} = \frac{c_2^{(1)}-\lambda_2T}{c_1^{(1)}+c_2^{(1)}-\lambda_2T}.
\]
By the same approach, we can also get the expression for $\Dstar$.
\Halmos
\endproof

Lemma \ref{lem:Dstar1} provides the closed form of the decision of Algorithm \Dtwo in the first period, $v_1^{(1)}$, and the optimal solution to \eqref{eq:optD1}, $\Dstar_1$. The following lemma shows their distance and the proof can be found in Appendix \ref{secappend:thm1}.

\begin{lemma} \label{lem:difference}
    $\mathbb{E}[v_1^{(1)}-\Dstar]=O(1/T)$.
\end{lemma}

Given that the initial value of the stochastic process $v_1^{(t)}$ is close to $\Dstar$, we discuss how $v_1^{(t)}$ evolves. This requires to split the analysis into two cases based on the remaining availability and time.

\text{ }

\subsubsection{Case 1: $c_1^{(t)}+c_2^{(t)} > (\lambda_1+\lambda_2)(T-t+1)$ }

\text{ }

In this case, by Lemma \ref{lem:Dstar1}, we can obtain
\[
v_1^{(t)}=\max\left\{\min\left\{\vstar_1,\frac{c_2^{(t)}}{\lambda_1(T-t+1)}-\frac{\lambda_2}{\lambda_1}\right\}, 1-\frac{c_1^{(t)}}{\lambda_1(T-t+1)}\right\}.
\]

We have $v_1^{(1)} = \max \left\{ \min\left\{ \vstar_1,\frac{c_2^{(1)}-\lambda_2 T}{\lambda_1 T}   \right\}, \frac{\lambda_1 T-c_1^{(1)}}{\lambda_1 T} \right\}$. If $\vstar_1$ is in the middle of the interval $\left[\frac{\lambda_1 T-c_1^{(1)}}{\lambda_1 T}, \frac{c_2^{(1)}-\lambda_2 T}{\lambda_1 T}\right]$, then $v_1^{(1)}=\vstar_1$. For each period $t$, $v_1^{(t)}$ is the projection of $\vstar_1$ into $\left[1-\frac{c_1^{(t)}}{\lambda_1(T-t+1)},  \frac{c_2^{(t)}}{\lambda_1(T-t+1)}-\frac{\lambda_2}{\lambda_1}\right]$. As the value of the lower and upper bound of the interval is slightly changing for each $t$, we have most likely, $\vstar_1$ is always in this interval, and $v_1^{(t)}=\vstar_1=v_1^{(1)}$. Therefore, in this case, we almost have zero loss.

However, if $\vstar_1$ is outside the interval $\left[\frac{\lambda_1 T-c_1^{(1)}}{\lambda_1 T}, \frac{c_2^{(1)}-\lambda_2 T}{\lambda_1 T}\right]$, for example, $\vstar_1 > \frac{c_2^{(1)}-\lambda_2 T}{\lambda_1 T}$, then $v_1^{(1)}=\frac{c_2^{(1)}-\lambda_2 T}{\lambda_1 T}$. Therefore, in each period $t$, if the upper bound moves left, we will have $v_1^{(t)}<v_1^{(1)}$, which result in some loss. The next lemma constructs a martingale to characterize the movement of the upper bound of the interval.

\begin{lemma}\label{lem:case1martingale}
    Let $\bigtriangleup_2(s)$ be the amount of type $2$ resources assigned in period $s$. Suppose that $v_1^{(1)}=\frac{c_2^{(1)}-\lambda_2 T}{\lambda_1 T}$. Construct a stochastic process $\epsilon^{U}_t=v_1^{(1)}-\left(\frac{c_2^{(t)}}{\lambda_1(T-t+1)}-\frac{\lambda_2}{\lambda_1}\right)$. Then, we have
    \[
    \epsilon^{U}_t = \frac{1}{\lambda_1}\sum_{s=1}^{t-1}\frac{\bigtriangleup_2(s)-\mathbb{E}[\bigtriangleup_2(s)]}{T-s+1},
    \]
    and $\epsilon^{U}_t$ is a martingale.
\end{lemma}

\proof{Proof of Lemma \ref{lem:case1martingale}}
As $v_1^{(1)}=\frac{c_2^{(1)}-\lambda_2 T}{\lambda_1 T}$, we have $c_2^{(1)}=\lambda_1Tv_1^{(1)}+\lambda_2T$.
If $t=2$, we have $\frac{c_2^{(1)}}{\lambda_1(T-2+1)}-\frac{\lambda_2}{\lambda_1}=\frac{c_2^{(1)}-\bigtriangleup_2(1)}{\lambda_1(T-1)}-\frac{\lambda_2}{\lambda_1}$. Then, we can obtain
\begin{align*}
    \epsilon^{U}_2 &= \frac{1}{\lambda_1}\left(\frac{c_2^{(1)}}{T}-\frac{c_2^{(1)}-\bigtriangleup_2(1)}{T-1} \right) = \frac{1}{\lambda_1}\frac{\bigtriangleup_2(1)T-c_2^{(1)}}{T(T-1)} \\&= \frac{1}{\lambda_1}\frac{\bigtriangleup_2(1)T-(\lambda_1Tv_1^{(1)}+\lambda_2T))}{T(T-1)} = \frac{1}{\lambda_1}\frac{\bigtriangleup_2(1)-\mathbb{E}[\bigtriangleup_2(1)]}{T-1},
\end{align*}
where the last step is because $\bigtriangleup_2(1)=1$ if and only if premium request arrives or basic request arrives and accepted the upgrading plan. Therefore, the expectation of $\bigtriangleup_2(1)$ is $\lambda_1v_1^{(1)}+\lambda_2$.

Next, by math induction, we have
\[
\epsilon^{U}_t=\frac{1}{\lambda_1}\sum_{s=1}^{t-1}\frac{\bigtriangleup_2(s)-\mathbb{E}[\bigtriangleup_2(s)]}{T-s+1}.
\]
Because at each $t$ expected increment $\mathbb{E}\left[ \bigtriangleup_2(t)-\mathbb{E}[\bigtriangleup_2(t)] \right] = 0$, we have that $\epsilon^{U}_t$ is a martingale.
\Halmos
\endproof
Similarly, if $\vstar_1 < \frac{c_2^{(1)}-\lambda_2 T}{\lambda_1 T}$, then $v_1^{(1)}=\frac{\lambda_1 T-c_1^{(1)}}{\lambda_1 T}$. The next lemma constructs a martingale to characterize the movement of $\frac{\lambda_1 (T-t+1)-c_1^{(t)}}{\lambda_1 (T-t+1)}$, and the proof can be found in Appendix \ref{secappend:thm1}.

\begin{lemma}\label{lem:case1martingaleL}
    Let $\bigtriangleup_1(s)$ be the amount of type $1$ resources assigned in period $s$. Suppose that $v_1^{(1)}=\frac{\lambda_1 T-c_1^{(1)}}{\lambda_1 T}$. Construct a stochastic process $\epsilon^{L}_t=v_1^{(1)}-\frac{\lambda_1 (T-t+1)-c_1^{(t)}}{\lambda_1 (T-t+1)}$. Then, we have
    \[
    \epsilon^{L}_t = \frac{1}{\lambda_1}\sum_{s=1}^{t-1}\frac{\mathbb{E}[\bigtriangleup_1(s)-\bigtriangleup_1(s)]}{T-s+1},
    \]
    and $\epsilon^{L}_t$ is a martingale.
\end{lemma}

Lemmas \ref{lem:case1martingale} and \ref{lem:case1martingaleL} show that in this case, $v_1^{(1)}-v_1^{(t)}$ is a martingale with expected increment $0$.

\subsubsection{Case 2: $c_1^{(t)}+c_2^{(t)} \leq (\lambda_1+\lambda_2)(T-t+1)$ }

\text{}

In this case, Lemma \ref{lem:Dstar1} implies that
\[
v_1^{(t)}=1-\frac{c_1^{(t)}}{c_1^{(t)}+c_2^{(t)}-\lambda_2 (T-t+1)}.
\]
However, different with case 1, it is easy to see that $v_1^{(1)}-v_1^{(t)}$ in this case is \textbf{not} a martingale. The next lemma constructs a complicated stochastic process $\alpha_t$ which satisfies two key properties: (i) it is a martingale with expected increment $0$, (ii) the difference between $\alpha_t$ and $\epsilon_t$ can be bounded by $O(1/T^2)$ for $t \in [1,\tau]$.

\begin{lemma} \label{lem:case2martingale}
    Define $\epsilon_t = v_1^{(1)}-v_1^{(t)}$ for any $t \in [T]$. 
    For any constant $\gamma \in (0,1)$, on the time horizon $[1, \gamma T]$, there exists a martingale $\alpha_t$ and a constant $\zeta$, such that $|\alpha_t - \epsilon_t| \leq \zeta / T^2$.
\end{lemma}

\proof{Proof of Lemma \ref{lem:case2martingale}}

First, when $c_1^{(t)}+c_2^{(t)} = (\lambda_1+\lambda_2)(T-t+1)$, we have
\[
1-\frac{c_1^{(t)}}{c_1^{(t)}+c_2^{(t)}-\lambda_2 (T-t+1)} = \frac{c_2^{(t)}}{\lambda_1(T-t+1)}-\frac{\lambda_2}{\lambda_1} = 1-\frac{c_1^{(t)}}{\lambda_1(T-t+1)}. 
\]
This implies that the phase transition during the decision-making process is smooth and continuous. By Lemmas \ref{lem:case1martingale} and \ref{lem:case1martingaleL}, we have when $c_1^{(t)}+c_2^{(t)} > (\lambda_1+\lambda_2)(T-t+1)$, $\epsilon_t$ is a martingale. Therefore, we only need to take $\alpha_t=\epsilon_t$. 

Then, we focus on the case where $v_1^{(t)}=1-\frac{c_1^{(t)}}{c_1^{(t)}+c_2^{(t)}-\lambda_2 (T-t+1)}$. We proceed by induction. Take $\alpha_1=0$. For $t=2$,
\begin{itemize}
    \item With probability $1-\lambda_1-\lambda_2$, no request arrives, and we have: 
    \begin{align*}
        v_1^{(1)}-v_1^{(2)}&=-\left(\frac{c_1^{(1)}}{c_1^{(1)}+c_2^{(1)}-\lambda_2 T}-\frac{c_1^{(1)}}{c_1^{(1)}+c_2^{(1)}-\lambda_2 T + \lambda_2}\right) \\&= -c_1^{(1)} \lambda_2 \frac{1}{(c_1^{(1)}+c_2^{(1)}-\lambda_2 T)(c_1^{(1)}+c_2^{(1)}-\lambda_2 T + \lambda_2)} \\&= -c_1^{(1)} \lambda_2 \frac{1}{(c_1^{(1)}+c_2^{(1)}-\lambda_2 T)^2} + O(1/T^3).
    \end{align*}
    \item With probability $\lambda_1(1-v_1^{(1)})$, a type $1$ request arrives and the upgrade is not accepted. Then:
    \begin{align*}
        v_1^{(1)}-v_1^{(2)} &=-\left(\frac{c_1^{(1)}}{c_1^{(1)}+c_2^{(1)}-\lambda_2 T}-\frac{c_1^{(1)}-1}{c_1^{(1)}-1+c_2^{(1)}-\lambda_2 T + \lambda_2}\right) \\
        &= - \frac{c_2^{(1)}-\lambda_2 T + \lambda_2 c_1^{(1)}}{(c_1^{(1)}+c_2^{(1)}-\lambda_2 T)(c_1^{(1)}-1+c_2^{(1)}-\lambda_2 T + \lambda_2)} \\
        &= -c_1^{(1)}\left(\frac{v_1^{(1)}}{1-v_1^{(1)}}+\lambda_2\right) \frac{1}{(c_1^{(1)}+c_2^{(1)}-\lambda_2 T)(c_1^{(1)}-1+c_2^{(1)}-\lambda_2 T + \lambda_2)} \\
        &= -c_1^{(1)}\left(\frac{v_1^{(1)}}{1-v_1^{(1)}}+\lambda_2\right)  \frac{1}{(c_1^{(1)}+c_2^{(1)}-\lambda_2 T)^2} + O(1/T^3),
    \end{align*}
    where the third equality is because
    \begin{align*}
        c_2^{(1)}-\lambda_2 T + \lambda_2 c_1^{(1)} &= \left(\frac{c_2^{(1)}-\lambda_2 T}{c_1^{(1)}}+\lambda_2 \right)c_1^{(1)} = (\frac{c_2^{(1)}-\lambda_2 T}{c_1^{(1)}+c_2^{(1)}-\lambda_2 T}\frac{c_1^{(1)}+c_2^{(1)}-\lambda_2 T}{c_1^{(1)}}+\lambda_2)c_1^{(1)} \\&= \left(\frac{v_1^{(1)}}{1-v_1^{(1)}}+\lambda_2\right)c_1^{(1)}.
    \end{align*}
    \item With probability $\lambda_2+\lambda_1v_1^{(1)}$, a type $2$ request arrives or a type $1$ request arrives and accepts the upgrading plan, we have
    \begin{align*}
        v_1^{(1)}-v_1^{(2)}&=\frac{c_1^{(1)}}{c_1^{(1)}+c_2^{(1)}-1-\lambda_2 T+\lambda_2}-\frac{c_1^{(1)}}{c_1^{(1)}+c_2^{(1)}-\lambda_2 T } \\&= c_1^{(1)} (1-\lambda_2) \frac{1}{(c_1^{(1)}+c_2^{(1)}-\lambda_2 T)(c_1^{(1)}+c_2^{(1)}-1-\lambda_2 T + \lambda_2)} \\&= c_1^{(1)} (1-\lambda_2)  \frac{1}{(c_1^{(1)}+c_2^{(1)}-\lambda_2 T)^2} + O(1/T^3).
    \end{align*}
\end{itemize}

Next, we define $\alpha_2$ as the value of $v_1^{(1)}-v_1^{(2)}$ without the $O(1/T^3)$ term:
\begin{align*}   \alpha_2 = \left\{ \begin{array}{ll}
         -c_1^{(1)} \lambda_2 \frac{1}{(c_1^{(1)}+c_2^{(1)}-\lambda_2 T)^2}  &\quad  \text{with probability } 1-\lambda_1-\lambda_2\\
        -c_1^{(1)}\left(\frac{v_1^{(1)}}{1-v_1^{(1)}}+\lambda_2\right)  \frac{1}{(c_1^{(1)}+c_2^{(1)}-\lambda_2 T)^2} &\quad  \text{with probability } \lambda_1(1-v_1^{(1)}) \\
        c_1^{(1)} (1-\lambda_2)  \frac{1}{(c_1^{(1)}+c_2^{(1)}-\lambda_2 T)^2} &\quad  \text{with probability } \lambda_2+\lambda_1v_1^{(1)} 
        \end{array} \right. \end{align*}

Then, we calculate the expectation of $\alpha_2$:
\begin{align*}
    \mathbb{E}[\alpha_2] = \frac{c_1^{(1)}}{(c_1^{(1)}+c_2^{(1)}-\lambda_2 T)^2}\left(-\lambda_2(1-\lambda_1-\lambda_2)-\lambda_1(1-v_1^{(1)}) + (\lambda_2+\lambda_1v_1^{(1)})(1-\lambda_2) \right) =0.
\end{align*}

Therefore, we define the stochastic process
\begin{align*}   \alpha_t = \left\{ \begin{array}{ll}
         -c_1^{(t)} \lambda_2 \frac{1}{(c_1^{(t)}+c_2^{(t)}-\lambda_2 (T-t+1))^2}  &\quad  \text{with probability } 1-\lambda_1-\lambda_2\\
        -c_1^{(t)}\left(\frac{v_1^{(t)}}{1-v_1^{(t)}}+\lambda_2\right)  \frac{1}{(c_1^{(t)}+c_2^{(t)}-\lambda_2 (T-t+1))^2} &\quad  \text{with probability } \lambda_1(1-v_1^{(t)}) \\
        c_1^{(t)} (1-\lambda_2)  \frac{1}{(c_1^{(t)}+c_2^{(t)}-\lambda_2 (T-t+1))^2} &\quad  \text{with probability } \lambda_2+\lambda_1v_1^{(t)} 
        \end{array} \right. \end{align*}

By math induction, we have $\alpha_t$ is a martingale with expectation $0$. In addition, between $\alpha_t$ and $\epsilon_t$, there is an error term of $O(\frac{1}{(T-t+1)^3})$. On time horizon $[1,\gamma T]$, we can use a union bound to upper bound the gap between $\alpha_t$ and $\epsilon_t$: for any  $t \in [1,\gamma T]$,  \[|\alpha_t - \epsilon_t| = \sum_{t=1}^{\gamma T}O\left(\frac{1}{(T-t+1)^3}\right) \leq \frac{\zeta}{T^2}\] for some constant $\zeta>0$.
\Halmos
\endproof

Lemma \ref{lem:case2martingale} shows that we can always find a martingale to approximate the gap between our decision variable $v_1^{(t)}$ and the optimal variable $v_1^{(1)}$ for any $t \in [1, \gamma T]$, where $\gamma \in (0,1)$. Recall that our goal is to approximate this gap for any $t \in [1,\tau]$. Therefore, we capture the value of the stopping time $\tau$ by the following lemma. 

\begin{lemma} \label{lem:stopping2}
For any arbitrarily small constant $\eta>0$, with probability at least $1-O(1/T)$,
    \[
    \left(\frac{c_1^{(1)}+c_2^{(1)}-\lambda_2 T}{\lambda_1 T}-\eta\right)T \leq \tau \leq \left(\frac{c_1^{(1)}+c_2^{(1)}-\lambda_2 T}{\lambda_1 T}+\eta\right)T
    \]
\end{lemma}

\proof{Proof of Lemma \ref{lem:stopping2}}
Let $\bigtriangleup_1(s)$ be the total number of type $1$ resources assigned in period $s$, and we have $\bigtriangleup_1(s) \sim \text{Ber}\left(\lambda_1(1-v_1^{(s)})\right)$. Then, we have
\begin{align*}
    \mathbb{P}\left(\tau \geq (\frac{c_1^{(1)}+c_2^{(1)}-\lambda_2 T}{\lambda_1 T}-\eta)T \right)&=\mathbb{P}\left(\sum_{s=1}^{(\frac{c_1^{(1)}+c_2^{(1)}-\lambda_2 T}{\lambda_1 T}-\eta)T}\bigtriangleup_1(s)<c_1^{(1)} \right) \\&= \mathbb{P}\left(\sum_{s=1}^{(\frac{c_1^{(1)}+c_2^{(1)}-\lambda_2 T}{\lambda_1 T}-\eta)T}\left(\bigtriangleup_1(s)-\lambda_1\frac{c_1^{(1)}}{c_1^{(1)}+c_2^{(1)}-\lambda_2 T} \right)< \eta T \right).
\end{align*}
 By Lemma \ref{lem:case2martingale}, we have $\mathbb{E}[\bigtriangleup_2(s)] = \lambda_1\frac{c_1^{(1)}}{c_1^{(1)}+c_2^{(1)}-\lambda_2 T} +O(1/T^2)$. Therefore, by Chernoff's inequality:
\[
\mathbb{P}\left(\tau \geq \left(\frac{c_1^{(1)}+c_2^{(1)}-\lambda_2 T}{\lambda_1 T}-\eta\right)T \right) < \frac{1}{T}.
\]
Similarly, with a symmetric statement, we can also get 
\[
\mathbb{P}\left(\tau \leq \left(\frac{c_1^{(1)}+c_2^{(1)}-\lambda_2 T}{\lambda_1 T}+\eta\right)T \right) < \frac{1}{T}.
\]
\Halmos
\endproof

Lemma \ref{lem:stopping2} shows that $\tau = \gamma T$ for some constant $\gamma \in (0,1)$. Therefore, by Lemma \ref{lem:case2martingale}, the error of the constructed martingale $\alpha_t$ can be bounded by $O(\frac{1}{T^2})$. 

Combined the results from \textit{Case 1} and \textit{Case 2}, we are ready to upper bound the value of

$\mathbb{E}\left[\sum_{t \in [\tau-1]}\left(R(\Dstar_1)-R(v_1^{(t)})\right) \right]$.

\begin{lemma}\label{lem:secondupper}
    \[
    \mathbb{E}\left[\sum_{t \in [\tau-1]}\left(R(\Dstar_1)-R(v_1^{(t)})\right) \right] = O(\log T).
    \]
\end{lemma}

The proof, detailed in Appendix \ref{secappend:thm1}, utilizes the Taylor series expansion of $R(v_1^{(t)})$. The underlying intuition draws from Lemmas \ref{lem:case1martingale}, \ref{lem:case1martingaleL}, and \ref{lem:case2martingale}, which establish the existence of a martingale closely approximating the stochastic process $v_1^{(t)}$. Given that $R(\cdot)$ is bounded, continuous, and concave, \cite{jasin2012re,jasin2014reoptimization} suggest that the expected difference $\mathbb{E}\left[\sum_{t \in [\tau-1]}\left(R(\Dstar_1)-R(v_1^{(t)})\right) \right]$ scales logarithmically with $T$, denoted as $O(\log T)$.

Finally, we bound $\mathbb{E}[W^{\textsc{U}}(\tau,T)-W^{\pi}(\tau,T)]$ by the proposition below.

\begin{proposition} \label{prop:tauregret}
The regret for the second time segment is: $\mathbb{E}[W^{\textsc{U}}(\tau,T)-W^{\pi}(\tau,T)] = O(1)$.
\end{proposition}

\proof{Proof of Proposition \ref{prop:tauregret}}
In the horizon $[\tau,T]$, we only have type $2$ resources. In this case, when a type $1$ request arrives, we can only upgrade them for free. Therefore, accepting a type $1$, $2$ request can generate a revenue of $r_1$, $r_2$ respectively. This turns the problem to be the classical quantity-based single leg revenue management problem under stochastic arrival process. By \citet{vera2021bayesian}, the regret is $O(1)$.
\Halmos
\endproof

Combining Lemma \ref{lem:secondupper} and Proposition \ref{prop:tauregret} in Equation \eqref{eq:regr2}, we can compute the regret.
\begin{align*}
    R(T) = O(\log T) + O(1) = O(\log T)%\\&\leq \mathbb{E}\left[\sum_{t=1}^{\tau-1}\left(R(\Dstar_1)-R(v_1^{(t)})\right) \right] + \mathbb{E}[W^{\textsc{U}}(\tau,T)-W^{\pi}(\tau,T)],
\end{align*}

\section{Algorithm for Multiple Types of Resources} \label{sec:multiple}

In this section, we focus on scenarios involving $n$ distinct types of resources, where $n$ is a given parameter. Similarly to the steps in the previous section, we first describe the upper bound of the performance of any online mechanism in Section \ref{subsec:mupper}, and subsequently, we present and analyze our algorithm in Section \ref{subsec:malg}.

\subsection{Upper Bound on Performance of Any Online Algorithm} \label{subsec:mupper}

In Section \ref{subsec:upper}, we introduced the hybrid programming \eqref{eq:upperdp}, whose optimal objective value establishes an upper bound across all algorithms. Nevertheless, the complexity of \eqref{eq:upperdp} arises from the necessity to discuss the sequence in which resources become exhausted. Considering the presence of $n$ distinct resource types in this context, it is impracticable to exhaustively enumerate every potential sequence of depletion for each resource type. To address this, we define $\mathbf{\Lambda}=[\Lambda_1,\Lambda_2,\ldots,\Lambda_n]$ as a vector of random variables, each representing the total number of arrivals for each resource type. Furthermore, $\mathbf{v}=[v_1,v_2,\ldots,v_{n-1}]$ is a vector of decision variables. Given $\Lambda_i$, $i \in [n]$, and a fixed upgrade probability $v_i$ for each type $i \in [n-1]$, we let the total revenue be $\WM(\mathbf{\Lambda},\mathbf{v})$. Therefore, the hybrid formulation can be articulated as:
\begin{align} \label{eq:optdefm}
    w^{\textsc{U}} \nonumber& = \max_{\mathbf{v}} \WM(\mathbf{\Lambda},\mathbf{v})    \tag{n-type HP}
\end{align}
Since $w^{\textsc{U}}$ is a random variable depending on the value of $\mathbf{\Lambda}$, we define $W^{\textsc{U}}=\mathbb{E}[w^{\textsc{U}}]$ as the expected hindsight total revenue. 
The hybrid programming formulation \eqref{eq:optdefm} finds the optimal fixed upgrading probabilities for each category, with the objective of maximizing the total revenue under $\Lambda_i$ arrivals of type $i$ requests. Specifically, for type $n$ resources, which are classified as the most superior resource, we should accept all such requests. This acceptance results in a residual availability of $(c_n-\Lambda_n)^{+}$ resources of type $n$, which are then exclusively allocated for potential upgrades of requesters classified under type $n-1$. Given the arrival of $\Lambda_{n-1}$ requesters of type $n-1$, the optimal upgrading probability, denoted as $v_{n-1}$, is computed utilizing Equation \eqref{eq:upperdp}. With a little bit abuse of notation, we denote $w^{\textsc{U}}$ as $w^{\textsc{U}}(c_1,c_2,\Lambda_1,\Lambda_2)$, explicitly indicating its dependence on the availability of resources and number of requesters. Consequently, the optimal decision variable $\Dstar_{n-1}$ is determined as $\Dstar_{n-1} = \text{argmax}_{v_{n-1}} w^{\textsc{U}}(c_{n-1},(c_n-\Lambda_n)^{+},\Lambda_{n-1},0)$. Subsequently, for type $n-2$ requests, the remaining resources of type $n-1$ are quantified as $\big(c_{n-1}-\frac{\Lambda_{n-1}}{1-\Dstar_{n-1}}\big)^{+}$, leading to the derivation of the optimal solution for $v_{n-2}$, which is formulated as $\Dstar_{n-2} = \text{argmax}_{v_{n-2}} w^{\textsc{U}}(c_{n-2},\big(c_{n-1}-\frac{\Lambda_{n-1}}{1-\Dstar_{n-1}}\big)^{+},\Lambda_{n-2},0)$. This methodology is iteratively applied across the spectrum of types, from type $n-1$ to type $1$. The subsequent theorem encapsulates the derived upper bound of this analysis where the proof can be found in Appendix \ref{append:multiple}.%, and Figure \ref{fig:nest} illustrates the nested structure inherent in this approach.
\begin{theorem} \label{thm:upperm}
    The optimal solution to \eqref{eq:optdefm} is: $\Dstar_{n-1} = \text{argmax}_{v_{n-1}} w^{\textsc{U}}(c_{n-1},(c_n-\Lambda_n )^{+},\Lambda_{n-1},0)$; and for $i \in \{n-2, n-3, \ldots, 1\}$, $\Dstar_{i} = \text{argmax}_{v_{i}} w^{\textsc{U}}(c_{i},\big(c_{i+1}-\frac{\Lambda_{i+1}}{1-\Dstar_{i+1}}\big)^{+},\Lambda_{i},0)$. Moreover, $W^{\textsc{U}}$ is an upper bound of the total revenue among all online algorithms, where $W^{\textsc{U}}=\mathbb{E}[w^{\textsc{U}}]$.
\end{theorem}
%\begin{figure}[!tb]
%\center
%\includegraphics[width=0.9\textwidth]{nest.jpg}
%\caption{Intuition to the optimal solution to \eqref{eq:optdefm}} 
%\label{fig:nest}
%\end{figure}

\subsection{Algorithm Description and Regret Analysis} \label{subsec:malg}

In this section, we present Algorithm \Dn, designed for $n$ resource types. By applying Theorem \ref{thm:upperm}, we can obtain the closed form of optimal solution vector to \eqref{eq:optdefm}: $\mathbf{\Dstar}=[\Dstar_1, \Dstar_2, \ldots, \Dstar_{n-1}]$. Next, for each type $i$ resource, we protect 
 $\big(c_{i+1}^{(1)}-\frac{\lambda_{i+1}T}{1-\Dstar_{i+1}} \big)^{+}$ for type $i-1$ requesters to upgrade. Thereafter, we implement Algorithm \ref{alg:price} to dynamically decide the upgrade pricing between each adjacent pair of request types. The next theorem summarizes the result.

\begin{algorithm}[!tb]
\caption{\Dn}
\label{alg:ntype}
\DontPrintSemicolon % To remove semicolons at the end of lines
\KwIn{time horizon $T$; initial availability vector $\mathbf{c}=[c_1^{(1)},c_2^{(1)},\ldots,c_n^{(1)}]$; arrival rate vector $\mathbf{\lambda}=[\lambda_1, \lambda_2, \ldots, \lambda_n]$; optimal solution vector to a defined optimization problem: $\mathbf{\Dstar}=[\Dstar_1, \Dstar_2, \ldots, \Dstar_{n-1}, 0]$.}
\For{$i \in \{1,2,\ldots,n-1 \}$}{
    Apply \Dtwo for type $i$ and type $i+1$ requests with input $T$, $c_i^{(1)}$, $\big(c_{i+1}^{(1)}-\frac{\lambda_{i+1}T}{1-\Dstar_{i+1}} \big)^{+}$, $\lambda_i$, $0$.\;
}
\end{algorithm}

\begin{theorem} \label{thm:ntype}
    The regret of Algorithm \Dn is bounded by $O(n\log T)$.
\end{theorem}

The proof can be found in Appendix \ref{append:multiple}. The total revenue, $ w^{\textsc{U}}$, can be decomposed into the sum of total revenue accrued between each consecutive pair of resources. Theorem \ref{thm:upperm} shows that the mechanism by which \Dn allocates protection to each category $i$ resource for the potential upgrade demands of category $i-1$ is the same as \eqref{eq:optdefm}. By Theorem \ref{thm:2type}, the regret between each adjacent pair of resource types is $O(\log T)$. This leads to an overall regret of $O(n \log T)$. %In practical applications, the number of category of upgradable resources typically remains a finite constant. Thus, when $n$ is finite, the regret incurred by Algorithm \Dn is still $O(\log T)$.

\section{Empirical Study for Hotel Upgrading} \label{sec:oracleexp}

In this section, we present an empirical study based on real data from the hospitality industry. %We apply Algorithm \Dn to this dataset and provide a comprehensive summary of the findings.\cpcomment{This sentence is not giving enough information so we might as well erase it. The reason is that we haven't introduced the data yet, the assumptions that we made, and any data pre-processing, and how our algorithm could be applied on it.}

\begin{figure}[!tb]
\center
\includegraphics[width=0.45\textwidth]{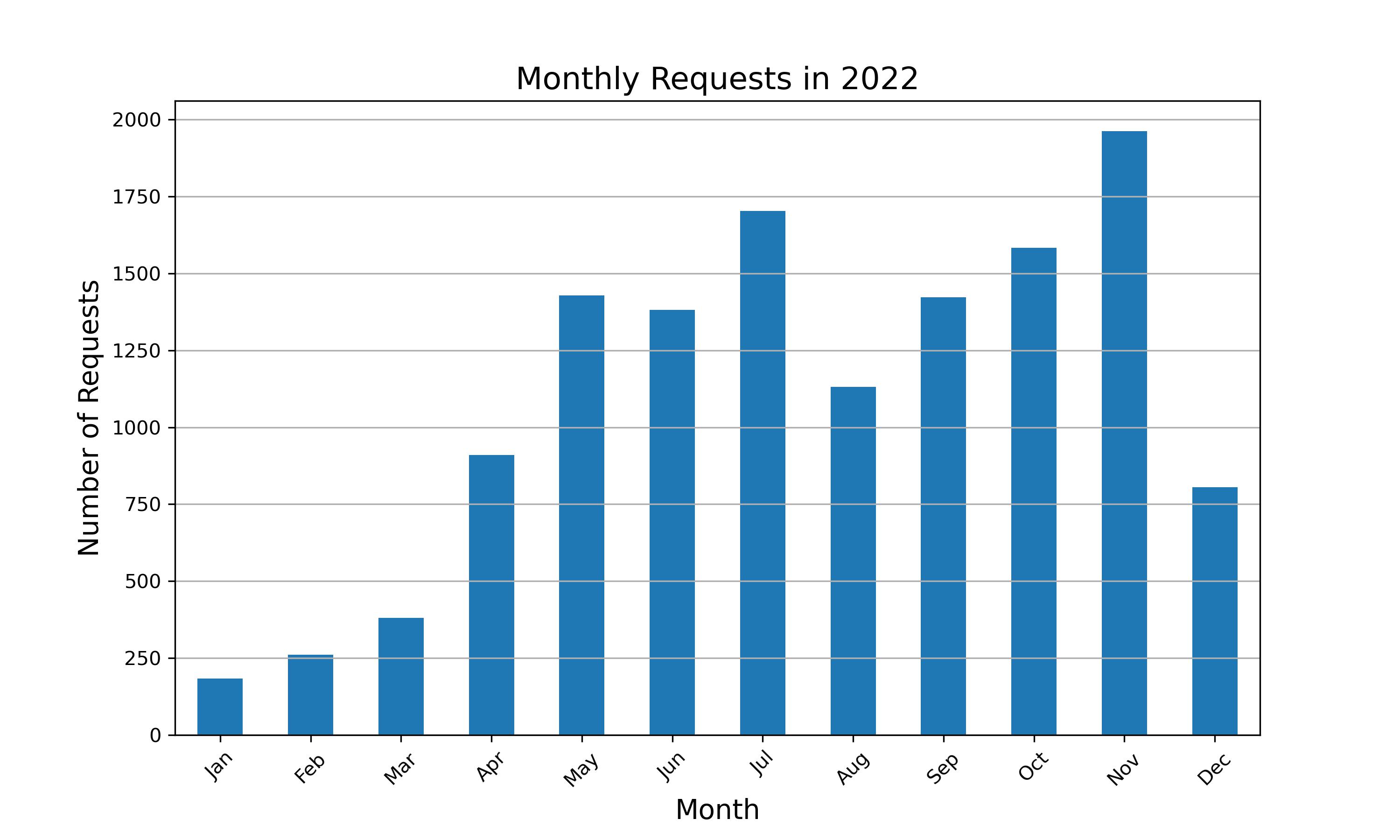}
\includegraphics[width=0.45\textwidth]{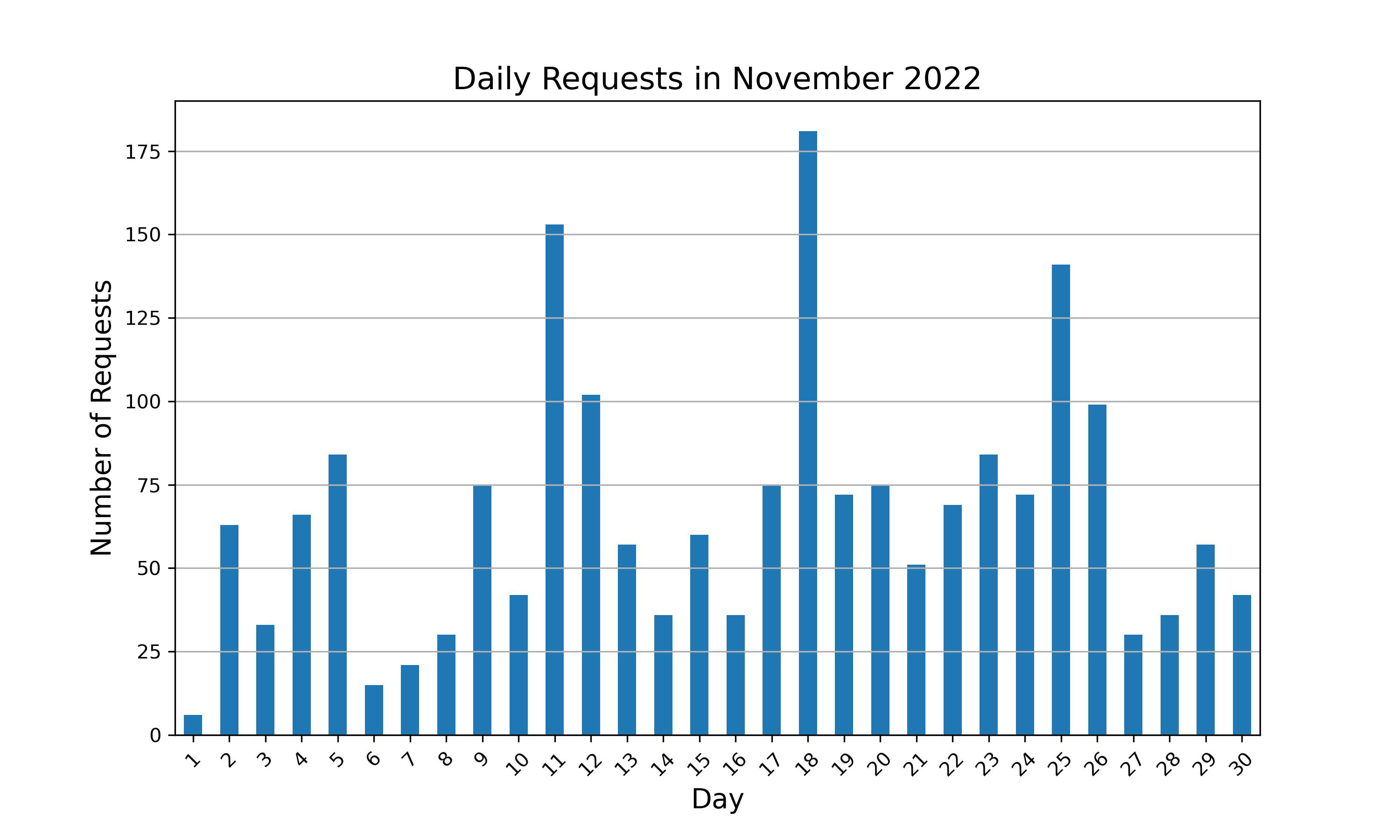}
\caption{Left: Distribution of number of requests in each month in 2022. Right: Distribution of a number of requests in each day in November 2022.} 
\label{fig:distribution}
\end{figure}

\xhdr{Data Description.} The dataset contains 13,155 upgrade offers recorded over one year. The left histogram of Figure \ref{fig:distribution} illustrates the monthly distribution of these offers. The peak in the volume of the offers is observed in November. Consequently, in this section, we use the data corresponding to November as an example for a detailed examination and implementation of the proposed algorithm\footnote{We also implement the algorithm for all days in all the months.}. The right histogram of Figure \ref{fig:distribution} describes the daily distribution of requests in November.

The dataset for each reservation request includes the following attributes: 
\begin{itemize}
    \item \textit{price}: initial payment by the customer.
    \item \textit{from}: requested room category.
    \item \textit{total}: cost of the next higher room category.
    \item \textit{to}: next higher room category.
    \item \textit{upgrade}: upgrade price presented to the customer.
    \item \textit{decision}: the customer accepts or rejects the upgrade offer.
\end{itemize}
It is important to note that room prices, even within the same category, are subject to dynamic fluctuation. To effectively measure the relative cost of an upgrade, the dataset employs the upgrade price proportion, calculated as $\textit{upgrade proportion}=\textit{upgrade}/(\textit{total}-\textit{price})$. In addition, the existing upgrading strategy does not incorporate information about the remaining capacities of each room category. Instead, it only considers factors related to pricing and demand. Therefore, for every request to room category $i$, the existing system provides roughly the same $\textit{upgrade proportion}$ for upgrading to the category $i+1$.

\xhdr{Simulation Setup.} The dataset encompasses three room categories with inventories of $60$, $30$, and $2$, respectively, within the examined booking channel. To derive functions $f_1(\cdot)$ and $f_2(\cdot)$ mapping the \textit{upgrade proportion} to the probability of accepting the upgrade plan for type 1 and type 2 requests, a regression analysis is conducted. Based on data visualization, we propose a parametric exponential decay family of functions, $\mathcal{F}$, consisting of all functions of the form $f(v)=e^{-av^{b}}$. Subsequent regression analysis identifies the functions with the minimal mean square error, resulting in $f_1(v_1)=e^{-4.4853 v_1^{0.9889}}$ and $f_2(v_2)=e^{-2.33 v_2}$, which best fit the observed data. 

We assume that each request corresponds to a single-night occupancy, which confines the decision-making process to a daily basis. Each day, we randomly generate $100$ permutations of the arrival order to represent arriving instances. The benchmark algorithm employs the upgrade price found in the data and follows a first-come, first-served protocol, continuing until either all requests have arrived or the inventory is full. To implement Algorithm \Dn, we first approximate the arrival rates, denoted as $\lambda_i$, $i \in [3]$. We set a horizon $T$ to simulate daily reservations maintained by hotels and add random noise to $\mathbf{\lambda}$, accounting for cancellation in guest attendance or arrivals without reservations: $\lambda_i \sim \text{Unif}(\Lambda_i-\sqrt{\Lambda_i}, \Lambda_i+\sqrt{\Lambda_i}) / T$. Subsequently, Algorithm \Dn is applied to the $100$ samples and we compare its performance to the benchmark using the same sample instances.

\begin{figure}[!tb]
\center
\includegraphics[width=0.45\textwidth]{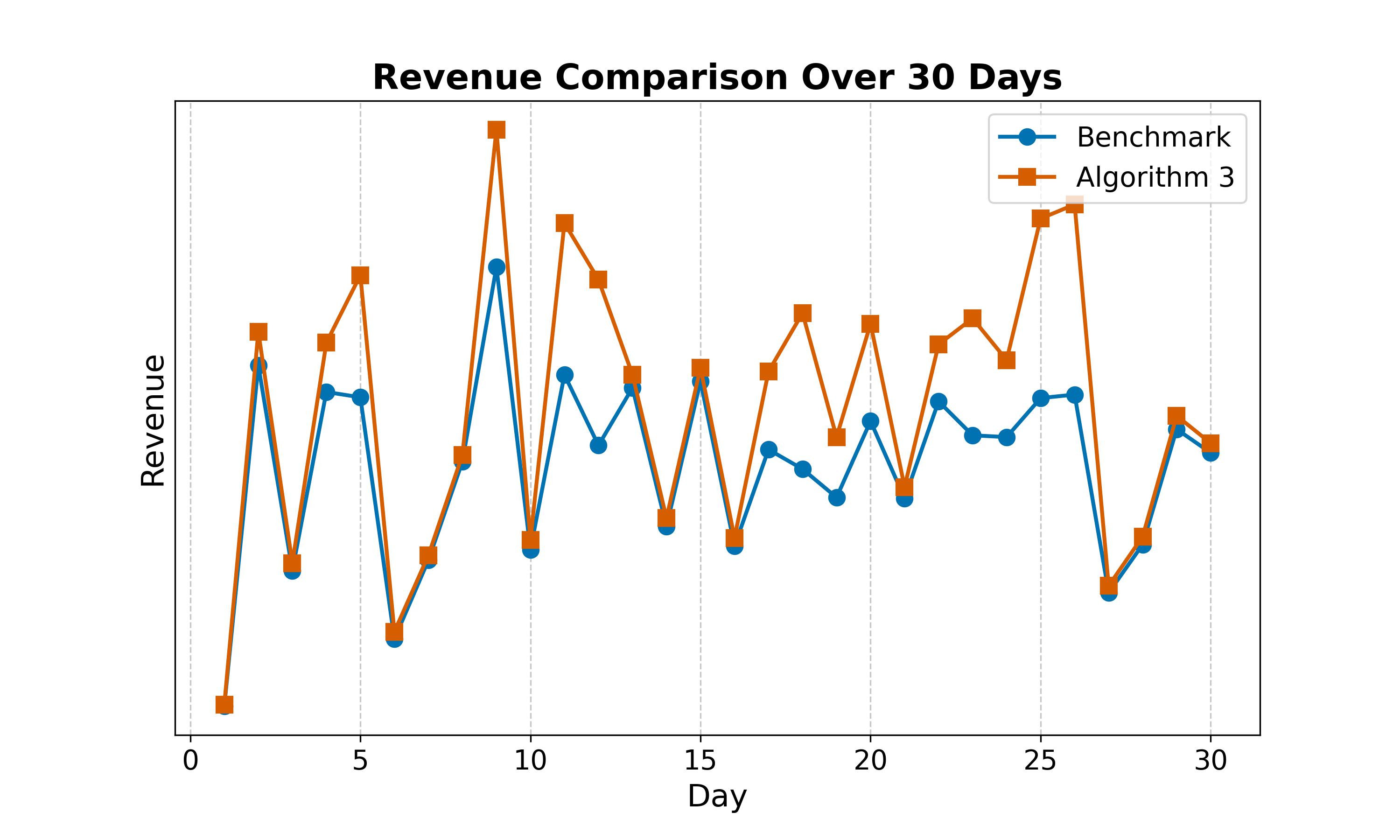}
\includegraphics[width=0.45\textwidth]{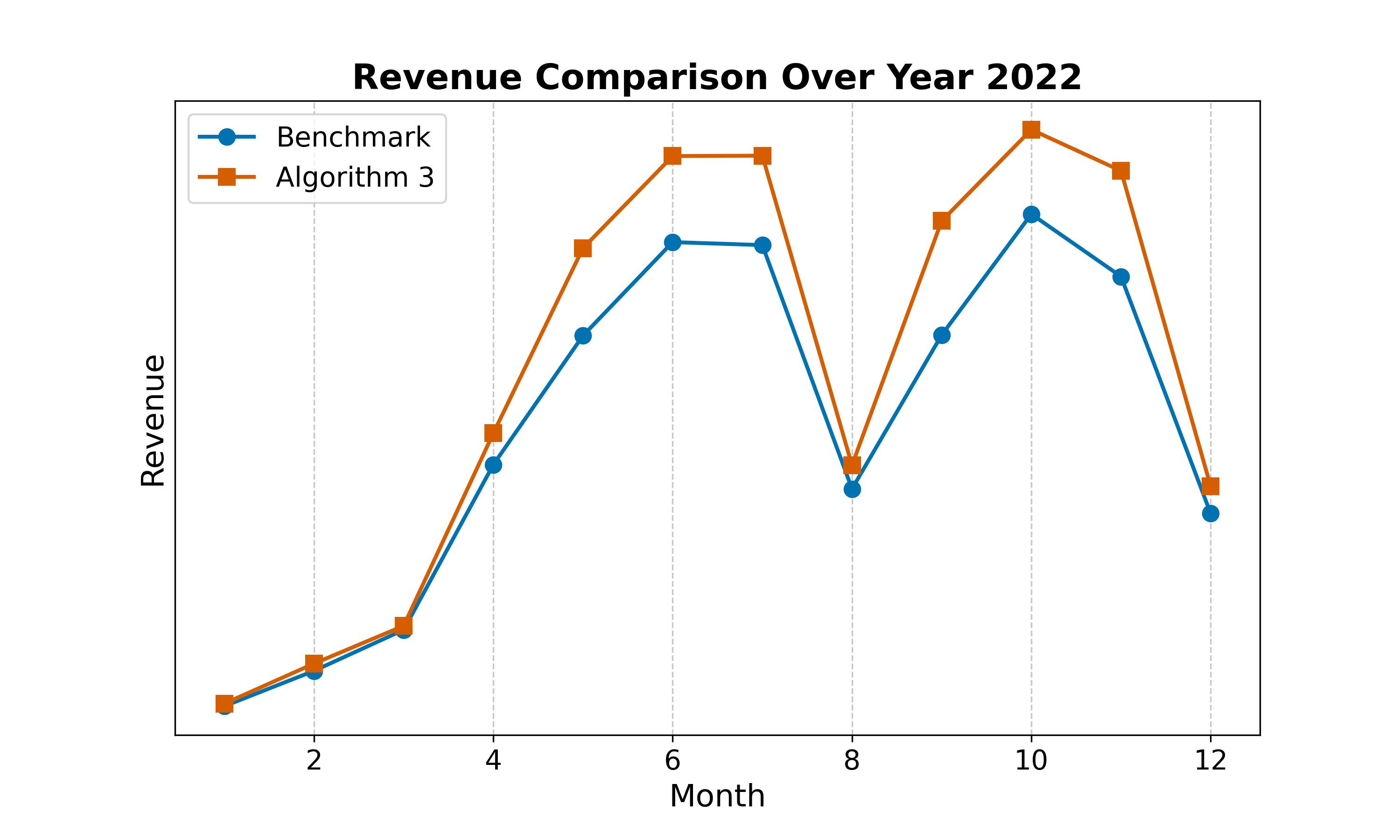}
\caption{Left: Daily total revenue of Benchmark and Algorithm \Dn in November 2022. Right: Monthly total revenue of Benchmark and Algorithm \Dn in 2022.} 
\label{fig:result2022}
\end{figure}

\xhdr{Results.} The left diagram of Figure \ref{fig:result2022} provides a comparative evaluation of daily revenues generated during November by both the benchmark and Algorithm \Dn. Despite noticeable fluctuations in daily demand, \Dn consistently outperforms the benchmark. Particularly on days characterized by higher demand, \Dn achieves a significant increase in total revenue, exceeding  $26\%$. The right panel of Figure \ref{fig:result2022} and Table \ref{table1} illustrate the monthly revenue trends throughout the year, with improvements of over $15\%$ highlighted in red. Over the entire period, the aggregate revenue generated by Algorithm \Dn exceeded the benchmark by approximately $17\%$. %\cp{I know from Google that the exact revenue values are usually proprietary information. Have we checked with Andrew that Oracle is fine with us reporting the exact dollar amounts? I also don't find the exact dollar amount of significant value in general. Maybe we can report the percentage change and then say that the actual revenue is between X and Y?{\color{blue}Jerry: As this is just a sample of data, and it's fine.}}

\iffalse
\begin{table}[ht]
\scriptsize
\caption{Monthly Summary in 2022}\label{table1}
\begin{center}
\begin{tabular}{ | m{2cm} | r | r | r |}
\hline
{\bf Month}   & {\bf Improvement}  \\
\hline
Jan  & 4.49\% \\
\hline
Feb  & 7.84\% \\
\hline
Mar  & 2.98\% \\
\hline
Apr  & 10.61\% \\
\hline
May  & {\color{red}20.56\%} \\
\hline
Jun  & {\color{red}16.75\%} \\
\hline
Jul  & {\color{red}17.24\%} \\
\hline
Aug  & 8.55\% \\
\hline
Sep  & {\color{red}27.22\%} \\
\hline
Oct  & {\color{red}15.20\%} \\
\hline
Nov  & {\color{red}26.82\%} \\
\hline
Dec  & 10.79\% \\
\hline
{\bf Total} & {\color{red}17.07\%} \\
\hline
\end{tabular}
\end{center}
\end{table}
\fi

\begin{table}[ht]
\caption{Monthly Summary in 2022}\label{table1}
\centering
\begin{tabular}{@{}l*{13}{r}@{}}
\toprule
\scriptsize \textbf{Month} & \scriptsize \textbf{Jan} & \scriptsize \textbf{Feb} & \scriptsize \textbf{Mar} & \scriptsize \textbf{Apr} & \scriptsize \textbf{\textcolor{red}{May}} & \scriptsize \textbf{\textcolor{red}{Jun}} & \scriptsize \textbf{\textcolor{red}{Jul}} & \scriptsize \textbf{Aug} & \scriptsize \textbf{\textcolor{red}{Sep}} & \scriptsize \textbf{\textcolor{red}{Oct}} & \scriptsize \textbf{\textcolor{red}{Nov}} & \scriptsize \textbf{Dec} & \scriptsize \textbf{\textcolor{red}{Total}} \\
\midrule
\scriptsize \textbf{Improvement} & \scriptsize 4.49\% & \scriptsize 7.84\% & \scriptsize 2.98\% & \scriptsize 10.61\% & \scriptsize \textcolor{red}{20.56\%} & \scriptsize \textcolor{red}{16.75\%} & \scriptsize \textcolor{red}{17.24\%} & \scriptsize 8.55\% & \scriptsize \textcolor{red}{27.22\%} & \scriptsize \textcolor{red}{15.20\%} & \scriptsize \textcolor{red}{26.82\%} & \scriptsize 10.79\% & \scriptsize \textcolor{red}{17.07\%} \\
\bottomrule
\end{tabular}
\end{table}

\section{Conclusion}

In this paper, we studied upgrading mechanisms for a fee for general resource allocation problems. Our proposed mechanism (Algorithm \Dn) achieves regret $O(\log T)$ against a ``hybrid programming'' benchmark that we analyze. Algorithm \Dn stands out not only for its theoretical regret bound but also for its performance in real-world hotel management scenarios.

There are a lot of avenues for future research. For instance, accurately estimating the functions $f_i(\cdot)$ presents a challenge in practical settings. This raises the question of whether effective blind online upgrading can be implemented without prior knowledge of $f_i(\cdot)$. Specifically, it necessitates developing strategies for upgrade pricing to make online learning of $f_i(\cdot)$ possible without incurring significant losses during the learning phase. Additionally, exploring other factors, such as the impact of customer loyalty on upgrading costs, is essential. Determining optimal strategies for discounting upgrade fees to reward loyalty requires careful consideration. Moreover, ensuring fairness in upgrading costs poses a significant challenge. Given that our dynamic upgrading mechanism assigns different fees to each request, incorporating fairness constraints necessitates adjustments to the algorithm's structure. Addressing these questions will enhance the efficacy of the online upgrading model and mechanism.  %We eagerly anticipate follow-up research that further refines and enhances the online upgrading mechanism.

\newpage
\bibliographystyle{ACM-Reference-Format}
\bibliography{reference}

\ECSwitch

\ECHead{Online Appendix}

\section{Supplemental Materials for Section \ref{subsec:upper}} \label{secappend:upper}

\begin{lemma} \label{lem:monotone}
    Fix any value of $\Lambda_1, \Lambda_2 \geq 0$, $\Lambda_2r_2 +\frac{c_2-\Lambda_2}{v_1} R_1(v_1) + r_1\min\Big\{\Lambda_1, c_1+c_2-\Lambda_2 \Big\}$ is a monotone decreasing function in $v_1 \in [0,1]$. $\Lambda_2r_2 + \frac{c_1}{1-v_1} R_1(v_1) + r_1\min\Big\{\Lambda_1, c_1+c_2-\Lambda_2 \Big\}$ is a monotone increasing function in $v_1 \in [0,1]$.
\end{lemma}

\proof{Proof of Lemma \ref{lem:monotone}}
First, to check the monotonicity of $\Lambda_2r_2 +\frac{c_2-\Lambda_2}{v_1} R_1(v_1) + r_1\min\Big\{\Lambda_1, c_1+c_2-\Lambda_2 \Big\}$, we only need to show that $\frac{R_1(v_1)}{v_1}$ is decreasing. This is obvious because $\frac{R_1(v_1)}{v_1}=p(v_1)$, which is defined as a monotone decreasing function.

Second, to check the monotonicity of $\Lambda_2r_2 + \frac{c_1}{1-v_1} R_1(v_1) + r_1\min\Big\{\Lambda_1, c_1+c_2-\Lambda_2 \Big\}$, we only need to show that $\frac{R_1(v_1)}{v_1}$ is increasing in $v_1$. As $R(\cdot)$ is concave, we have $(v_1 p(v_1))^{''}<0$, which implies that $p^{''}(v_1)<\frac{-2p^{'}(v_1)}{v_1}$. Then, we have
\begin{align*}
    \left(\frac{R(v_1)}{1-v_1}\right)^{'} &= \frac{p(v_1)+v_1(1-v_1)p^{'}(v_1)}{(1-v_1)^2}.
\end{align*}

Observe that
\[
\left(p(v_1)+v_1(1-v_1)p^{'}(v_1) \right)^{'} = 2(1-v_1)p^{'}(v_1)+p^{''}(v_1)(v_1-v_1^2)<2(1-v_1)p^{'}(v_1)+\frac{-2p^{'}(v_1)}{v_1}(v_1-v_1^2)<0,
\]
which implies that $p(v_1)+v_1(1-v_1)p^{'}(v_1)$ is a decreasing function. Therefore, 
\begin{align*}
    \left(\frac{R(v_1)}{1-v_1}\right)^{'} &= \frac{p(v_1)+v_1(1-v_1)p^{'}(v_1)}{(1-v_1)^2} \\&\geq \frac{p(1)+1 \cdot (1-1)p^{'}(1)}{(1-v_1)^2} = 0,
\end{align*}
and we can obtain that $\frac{R(v_1)}{1-v_1}$ is increasing. 
\Halmos
\endproof

\subsection{Proof of Theorem \ref{thm:sameopt}}

\proof{Proof of Theorem \ref{thm:sameopt}}
We split the proof into two parts, where in part 1, we show that the optimal solution to \ref{eq:upperdp} is the same as the one to \ref{eq:optD1}, and in part 2, we prove that $W^{\textsc{U}} \geq \mathbb{E}[W^{\textsc{opt}}]$.

\underline{\textit{Part 1}}

Let $\Dstar_1 = \text{argmax}_{v_1} \Lambda_2r_2 + \min\Big\{\Lambda_1, \frac{c_2-\Lambda_2}{v_1}, \frac{c_1}{1-v_1} \Big\}\Big( R_1(v_1)+r_1 \Big)$. Then, we show that $\Dstar_1$ is also the optimal solution to $w^{\textsc{U}}$. 

Define $\vstar_1(\textsc{U}_1) = \text{argmax}_{v_1 \in \mathcal V_1} w^{\textsc{U}_1}$, $\vstar_1(\textsc{U}_2) = \text{argmax}_{v_1 \in \mathcal V_2} w^{\textsc{U}_2}$, and $\vstar_1(\textsc{U}_3)=\text{argmax}_{v_1 \in \mathcal V_3} w^{\textsc{U}_3}$. By Lemma \ref{lem:monotone}, we have $\vstar_1(\textsc{U}_2)=\inf \mathcal V_2$ and $\vstar_1(\textsc{U}_3)=\sup \mathcal V_3$. Then, we split the discussion into two cases.

\textit{Case 1: $\mathcal V_1 = \emptyset$. } In this case, as $\frac{c_2-\Lambda_2}{v_1}$ is decreasing in $v_1$ and $\frac{c_1}{1-v_1}$ is increasing in $v_1$, these two functions can have only one intersection. Due to $\mathcal V_1 = \emptyset$, we have $\inf \mathcal V_2= \sup \mathcal V_3$, which implies that $\vstar_1(\textsc{U}_2)=\vstar_1(\textsc{U}_3)$. By definition, we have $w^{\textsc{U}_2} = w^{\textsc{U}_3} = w^{\textsc{U}}$. 

Moreover, by Lemma \ref{lem:monotone}, if $\mathcal V_1 = \emptyset$, we have $\Dstar_1 = \inf \mathcal V_2= \sup \mathcal V_3 = \vstar_1(\textsc{U}_2)=\vstar_1(\textsc{U}_3)$. Therefore, both $w^{\textsc{U}}$ and $W^{\textsc{h}}$ share the same optimal solution in this case.

\textit{Case 2: $\mathcal V_1 \neq \emptyset$. } In this case, as $\Lambda_1$ is a constant function, it can have only one intersection with $\frac{c_2-\Lambda_2}{v_1}$ and one intersection with $\frac{c_1}{1-v_1}$. Therefore, we have $\vstar_1(\textsc{U}_2) = \inf \mathcal V_2= \sup \mathcal V_1$ and $\vstar_1(\textsc{U}_3)=\sup \mathcal V_3$. This implies that 
\[
    w^{\textsc{U}_2} = \Lambda_2r_2 +\Lambda_1 R_1(\vstar_1(\textsc{U}_2)) + r_1\min\Big\{\Lambda_1, c_1+c_2-\Lambda_2 \Big\},
\]
and
\[
    w^{\textsc{U}_3} = \Lambda_2r_2 +\Lambda_1 R_1(\vstar_1(\textsc{U}_3)) + r_1\min\Big\{\Lambda_1, c_1+c_2-\Lambda_2 \Big\}.
\]
It is easy to see that 
\[
        w^{\textsc{U}_1} = \max_{v_1 \in \mathcal V_1} \Lambda_2r_2 + \Lambda_1 \Big( R_1(v_1)+r_1 \Big) \geq w^{\textsc{U}_2},
\]
and
\[
        w^{\textsc{U}_1} = \max_{v_1 \in \mathcal V_1} \Lambda_2r_2 + \Lambda_1 \Big( R_1(v_1)+r_1 \Big) \geq w^{\textsc{U}_3}.
\]
Therefore, we have $w^{\textsc{U}}=w^{\textsc{U}_1}$, and its optimal solution is $\vstar_1(\textsc{U}_1)$.

Moreover, by Lemma \ref{lem:monotone}, if $\mathcal V_1 \neq \emptyset$, we have $\Lambda_1 = \min\Big\{\Lambda_1, \frac{c_2-\Lambda_2}{\Dstar_1}, \frac{c_1}{1-\Dstar_1} \Big\}$, and this implies that $\Dstar_1 \in \mathcal V_1$. Therefore, $\Dstar_1 = \vstar_1(\textsc{U}_1)$, and both $w^{\textsc{U}}$ and $W^{\textsc{h}}$ share the same optimal solution in this case.

\underline{\textit{Part 2}}

Let $\Dstar_1$ be the optimal solution to  \eqref{eq:optD1}. From part 1, $\Dstar_1$ is also the optimal solution to \eqref{eq:upperdp}. Then, we split the proof into three cases:

\textit{Case 1: $\Lambda_1= \min\Big\{\Lambda_1, \frac{c_2-\Lambda_2}{\Dstar_1}, \frac{c_1}{1-\Dstar_1} \Big\}$. } In this case, by definition, $w^{\textsc{U}} = W^{\textsc{h}}$. In expectation, no resource is depleted during the time horizon $[0,T]$. As our decision process is Markovian and the objective is concave, by a well-known result (\cite{gallego1994optimal} Theorem II), we have $w^{\textsc{U}} = W^{\textsc{h}} \geq W^{\textsc{opt}}$.

\textit{Case 2: $\frac{c_2-\Lambda_2}{\Dstar_1}= \min\Big\{\Lambda_1, \frac{c_2-\Lambda_2}{\Dstar_1}, \frac{c_1}{1-\Dstar_1} \Big\}$. } In this case, the premium resource (type $2$) will be depleted at certain point during the time horizon $[0,T]$. Again, by \cite{gallego1994optimal} Theorem II, if the decision process stops when the premium resource is depleted, $\max_{v_1 \in \mathcal V_2} \Lambda_2r_2 +\frac{c_2-\Lambda_2}{v_1} (R_1(v_1)+r_1)$ serves an upper bound to the optimal policy. However, for a given 
$v_1 \in \mathcal V_2$, there is $c_1-\frac{c_2-\Lambda_2}{v_1}(1-v_1)$ basic resource (type $1$) remaining unused. Also, there is $\Lambda_1 - \frac{c_2-\Lambda_2}{v_1}$ type $1$ demand units who will arrive in the remaining time horizon in expectation. Therefore, the upper bound formulation is:
\begin{align*}
    & \text{ } \text{ } \text{ } \max_{v_1 \in \mathcal V_2} \Lambda_2r_2 +\frac{c_2-\Lambda_2}{v_1} (R_1(v_1)+r_1) + r_1 \min\Big\{\Lambda_1 - \frac{c_2-\Lambda_2}{v_1},  c_1-\frac{c_2-\Lambda_2}{v_1}(1-v_1) \Big\} \\& = \max_{v_1 \in \mathcal V_2} \Lambda_2r_2 +\frac{c_2-\Lambda_2}{v_1} R_1(v_1) + r_1\min\Big\{\Lambda_1, c_1+c_2-\Lambda_2 \Big\} \\& =w^{\textsc{U}_2} 
\end{align*}

Because $w^{\textsc{U}} \geq w^{\textsc{U}_2}$, we have $w^{\textsc{U}} \geq W^{\textsc{opt}}$.

\textit{Case 3: $\frac{c_1}{1-\Dstar_1}= \min\Big\{\Lambda_1, \frac{c_2-\Lambda_2}{\Dstar_1}, \frac{c_1}{1-\Dstar_1} \Big\}$. } In this case, the basic resource (type $1$) will be depleted at certain point during the time horizon $[0,T]$. Again, by \cite{gallego1994optimal} Theorem II, if the decision process stops when the basic resource is depleted, $\max_{v_1 \in \mathcal V_3} \Lambda_2r_2 +\frac{c_1}{1-v_1} (R_1(v_1)+r_1)$ serves an upper bound to the optimal policy. However, for a given 
$v_1 \in \mathcal V_3$, there is $c_2- \Lambda_2 - \frac{c_1}{1-v_1} v_1$ premium resource (type $2$) remaining unused. Also, there is $\Lambda_1 - \frac{c_1}{1-v_1}$ type $1$ demand units who will arrive in the remaining time horizon in expectation. These type $1$ demand units can be upgraded to the premium resource even if the basic resource is depleted. Thus, the upper bound formulation is:
\begin{align*}
    & \text{ } \text{ } \text{ } \max_{v_1 \in \mathcal V_3} \Lambda_2r_2 +\frac{c_1}{1-v_1} (R_1(v_1)+r_1) + r_1 \min\Big\{\Lambda_1 - \frac{c_1}{1-v_1}, c_2- \Lambda_2 - \frac{c_1}{1-v_1} v_1 \Big\} \\& = \max_{v_1 \in \mathcal V_3} \Lambda_2r_2 +\frac{c_1}{1-v_1} R_1(v_1) + r_1\min\Big\{\Lambda_1, c_1+c_2-\Lambda_2 \Big\} \\& =w^{\textsc{U}_3} 
\end{align*}

Because $w^{\textsc{U}} \geq w^{\textsc{U}_3}$, we have $w^{\textsc{U}} \geq W^{\textsc{opt}}$. This also implies that $W^{\textsc{U}} \geq \mathbb{E}[W^{\textsc{opt}}]$.
\Halmos
\endproof

\section{Supplementary Materials for Section \ref{subsec:alg2}} \label{secappend:thm1}

\proof{Proof of Lemma \ref{lem:firstupper}}
    Denote $\mathcal{I}_1$ as the time index that a type $1$ customer arrives before $\tau$.
\begin{align*}
     &\mathbb{E}[W^{\textsc{U}}(1,\tau-1)-W^{\pi}(1,\tau-1)] \\&= \mathbb{E}\left[\Lambda_2(1,\tau-1)r_2+\Lambda_1(1,\tau-1)\left(R_1(\Dstar_1)+r_1 \right)-\Lambda_2(1,\tau-1)r_2+\sum_{t \in \mathcal{I}_1}\left(R(v_1^{(t)})+r_1\right)  \right] \\&\leq \mathbb{E}\left[\sum_{t=1}^{\tau-1}\left(R(\Dstar_1)-R(v_1^{(t)})\right) \right],
\end{align*}
where the first equality is because by Wald's equality, at stopping time $\tau$, the remaining capacity of type $2$ resource is $c_2^{(\tau)} \approx \lambda_2(T-\tau+1)>0$, which implies that before $\tau$, no resource is depleted.     
\endproof

\proof{Proof of Lemma \ref{lem:difference}}
    If both $v_1^{(1)}$ and $\Dstar$ are $\vstar_1$, we have $v_1^{(1)}-\Dstar=0$. If $v_1^{(1)}=\frac{c_2^{(1)}-\lambda_2T}{\lambda_1 T}$ and $\Dstar = \frac{c_2^{(1)}-\Lambda_2}{\Lambda_1}$, since $\Lambda_1$ and $\Lambda_2$ are independent, we have
    \[
    \mathbb{E}[v_1^{(1)}-\Dstar]=\mathbb{E}\left[\frac{c_2^{(1)}-\lambda_2T}{\lambda_1 T}-\frac{c_2^{(1)}-\Lambda_2}{\Lambda_1}\right] = (c_2^{(1)}-\lambda_2T) \left( \frac{1}{\lambda_1 T} - \mathbb{E}\left[\frac{1}{\Lambda_1} \right] \right).
    \]
    As $\Lambda_1 \sim \text{Bin}(T,\lambda_1)$, we have
    \[
    \mathbb{E}\left[\frac{1}{\Lambda_1} \right] \approx \mathbb{E}\left[\frac{1}{1+\Lambda_1} \right].
    \]
    Use the fact that $\mathbb{E}\left[\frac{1}{a+X} \right] = \int_{0}^{1}t^{a-1}P_X(t)dt$, where $P_X(t)$ is the probability generating function for $X$, we can obtain
    \begin{align*}
        \mathbb{E}\left[\frac{1}{1+\Lambda_1} \right] &= \int_{0}^{1}t^{0}P_{\Lambda_1}(t)dt \\&= \int_{0}^{1} (1-\lambda_1+\lambda_1t)^T dt \\&= \frac{1-(1-\lambda_1)^{T+1}}{(T+1)\lambda_1}
    \end{align*}
    Therefore, we have $\frac{1}{\lambda_1 T} - \mathbb{E}\left[\frac{1}{\Lambda_1} \right] \approx \frac{1}{\lambda_1 T} -  \frac{1-(1-\lambda_1)^{T+1}}{(T+1)\lambda_1} = O(T^{-2})$. Similarly, if $v_1^{(1)}=\frac{\lambda_1 T-c_1^{(1)}}{\lambda_1 T}$ and $\Dstar_1=\frac{\Lambda_1 -c_1^{(1)}}{\Lambda_1}$, we can use the same method to derive that $\mathbb{E}[v_1^{(1)}-\Dstar]=O(T^{-2})$.

    Finally, if $v_1^{(1)} = \frac{c_2^{(1)}-\lambda_2T}{c_1^{(1)}+c_2^{(1)}-\lambda_2T}$ and $\Dstar_1=\frac{c_2^{(1)}-\Lambda_2}{c_1^{(1)}+c_2^{(1)}-\Lambda_2}$, we have
    \begin{align*}
    \mathbb{E}[v_1^{(1)}-\Dstar]&=\mathbb{E}\left[1-\frac{c_1^{(1)}}{c_1^{(1)}+c_2^{(1)}-\lambda_2T}-\left(1-\frac{c_1^{(1)}}{c_1^{(1)}+c_2^{(1)}-\lambda_2T}\right)\right] \\&= c_1^{(1)}\mathbb{E}\left[ \frac{1}{c_1^{(1)}+c_2^{(1)}-\Lambda_2} - \frac{1}{c_1^{(1)}+c_2^{(1)}-\lambda_2T} \right].
    \end{align*}
    As $\Lambda_2 \sim \text{Bin}(T,\lambda_2)$, by Hoeffding's inequality, with probability at least $1-1/T$, we have $\Lambda_2 \in [\lambda_2T-\sqrt{T}\log T, \lambda_2T+\sqrt{T}\log T ]$. Therefore, we have
    \begin{align*}
        \mathbb{E}\left[ \frac{1}{c_1^{(1)}+c_2^{(1)}-\Lambda_2} - \frac{1}{c_1^{(1)}+c_2^{(1)}-\lambda_2T} \right] &= (1-\frac{1}{T}) O(T^{-\frac{3}{2}}\log T) + \frac{1}{T} = O(\frac{1}{T}). 
    \end{align*}

\endproof

\proof{Proof of Lemma \ref{lem:case1martingaleL}}
As $v_1^{(1)}=\frac{\lambda_1 T-c_1^{(1)}}{\lambda_1 T}$, we have $c_1^{(1)}=\lambda_1T(1-v_1^{(1)})$.
If $t=2$, we have $\frac{\lambda_1 (T-2+1)-c_1^{(2)}}{\lambda_1 (T-2+1)}=1-\frac{c_1^{(1)}-\bigtriangleup_1(1)}{\lambda_1(T-1)}$. Then, we can obtain
\begin{align*}
    \epsilon^{L}_2 &= \frac{1}{\lambda_1}\left(\frac{c_1^{(1)}-\bigtriangleup_1(1)}{T-1}-\frac{c_1^{(1)}}{T} \right) = \frac{1}{\lambda_1}\frac{\lambda_1 T (1-v_1^{(1)})-\bigtriangleup_1(1) T}{T(T-1)} \\&=  \frac{1}{\lambda_1}\frac{\mathbb{E}[\bigtriangleup_1(1)]-\bigtriangleup_1(1)}{T-1},
\end{align*}
where the last step is because $\bigtriangleup_1(1)=1$ if and only if basic request arrives and rejected the upgrading plan. Therefore, the expectation of $\bigtriangleup_1(1)$ is $\lambda_1(1-v_1^{(1)})$.

Next, by math induction, we have
\[
\epsilon^{L}_t = \frac{1}{\lambda_1}\sum_{s=1}^{t-1}\frac{\mathbb{E}[\bigtriangleup_1(s)]-\bigtriangleup_1(s)}{T-s+1}.
\]
Because each time the expected increment $\mathbb{E}\left[\mathbb{E}[\bigtriangleup_1(t)]- \bigtriangleup_1(t) \right] = 0$, we immediately have $\epsilon^{L}_t$ is a martingale.
\Halmos
\endproof

\proof{Proof of Lemma \ref{lem:secondupper}}
Take $\zeta=\frac{1}{T^3}$, and we perturb the term $R(\Dstar_1)$ by $\zeta$ units to make $R$ differentiable at $\Dstar_1-\zeta$. In Lemma \ref{lem:difference}, we upper bound the expected difference of $\beta = \Dstar_1-v_1^{(1)}$ by $O(\frac{1}{T})$. Moreover, in Lemma \ref{lem:case2martingale}, we defined a martingale $\alpha_t$ which approximately captures $\epsilon_t=v_1^{(1)}-v_1^{(t)}$, and $|\alpha_t-\epsilon_t|=O(\frac{1}{T^2})$. By Taylor's expansion, we have
\begin{align*}
\sum_{t=1}^{\tau-1}R(v_1^{(t)}) &\geq (\tau-1)R(\Dstar_1-\zeta) + R'(\Dstar_1-\zeta)\sum_{t=2}^{\tau-1}(\beta+\epsilon_t-\zeta)+\frac{1}{2} R''(\Dstar_1-\zeta)\sum_{t=2}^{\tau-1}(\beta+\epsilon_t-\zeta)^2\\&+\sum_{t=1}^{\tau - 1}(R(v_1^{(t)})-\mathbb{E}[R(v_1^{(t)}) | \mathcal{F}_t]).
\end{align*}

Observe that $\sum_{s=1}^{t}(R(v_1^{(s)})-\mathbb{E}[R(v_1^{(s)}) | \mathcal{F}_s])$ is a martingale. By stopping time theorem, we have 
\[
\mathbb{E} \left[\sum_{t=1}^{\tau - 1}(R(v_1^{(t)})-\mathbb{E}[R(v_1^{(t)}) | \mathcal{F}_t]) \right]=0.
\]
By continuity of $R$, we have 
\[\mathbb{E}\left[(\tau-1) R(\Dstar_1) -  (\tau-1) R(\Dstar_1-\zeta)  \right]=O(\frac{1}{T^2}).
\]
As $\epsilon_t$ is close to a martingale $\alpha_t$, and $|\alpha_t-\epsilon_t|=O(\frac{1}{T^2})$, we have $\mathbb{E}[\sum_{t=2}^{\tau-1}\epsilon_t]=O(\frac{1}{T})$ and by Lemma \ref{lem:difference}, $\mathbb{E}[\beta]=O(\frac{1}{T})$, we can obtain
\[
-\mathbb{E}\left[ R'(\Dstar_1-\zeta)\sum_{t=2}^{\tau-1}(\beta+\epsilon_t-\zeta)\right] = O(1).
\]
Finally, we also have
\begin{align*}
    -\mathbb{E}[\frac{1}{2} R''(\Dstar_1-\zeta)\sum_{t=2}^{\tau-1}(\beta+\epsilon_t-\zeta)^2]  &\leq O(1)+ \frac{1}{2}R''(\Dstar_1-\zeta)\mathbb{E}\left[\sum_{t=2}^{\tau-1} (\epsilon_t)^2 \right] \\&=O(1)+ \frac{1}{2}R''(\Dstar_1-\zeta)\mathbb{E}[\sum_{t=2}^{\tau-1}\sum_{1 \leq s,v \leq t} \frac{1}{(T-s)(T-v)} \epsilon_s \epsilon_v  ]  \\&= O(1)+ \frac{1}{2}R''(\Dstar_1-\zeta)\mathbb{E}[\sum_{t=2}^{\tau-1}\sum_{s=1}^{t-1} \frac{1}{(T-s)^2} (\epsilon_s)^2  ] \\&= O(\log T).
\end{align*}
\Halmos
\endproof

\section{Supplementary Materials for Section \ref{sec:multiple}} \label{append:multiple}

\proof{Proof of Theorem \ref{thm:upperm}}
Firstly, the social cost of type $i+1$ dominates the maximum upgrade social cost of type $i$, namely $r_{i+1} \geq r_i+u_i$, we can decompose the hindsight formulation with multiple types to $n$ hindsight formulation with $1$ type of demand and $2$ types of resources. Therefore, by Theorem \ref{thm:sameopt}, the hybrid programming $W^{\textsc{U}}$ serves as an upper bound.

Secondly, to show the closed form solution of \eqref{eq:optdefm}, let us focus initially on type $n$ requests. As there are $\Lambda_n$ such requests and considering that $r_n > r_{n-1} > \ldots > r_1$, the optimal approach entails accepting $\min\{\Lambda_n ,c_n\}$ type $n$ requests. This results in a surplus of $(c_n-\Lambda_n )^{+}$ units of type $n$ resources. These remaining resources are exclusively allocated for fulfilling upgraded type $n-1$ demands. Given that each type $n-1$ demand utilizes either a unit of type $n-1$ or type $n$ resource, and in both cases, the total revenue surpasses that of other resources, it is necessary to accommodate as many type $n-1$ demands as feasible. Thus, by Theorem \ref{thm:sameopt}, the optimal upgrade probability of type $n-1$ demands is $\Dstar_{n-1} = \text{argmax}_{v_{n-1}}W^{\textsc{U}}(c_{n-1},(c_n-\Lambda_n )^{+},\Lambda_{n-1},0)$.

Next, by setting the upgrade probability of type $n-1$ demand as $\Dstar_{n-1}$, the surplus of type $n-1$ resources is $\left(c_{n-1}-\frac{\Lambda_{n-1}}{1-\Dstar_{n-1}}\right)^{+}$, and by similar reasoning, the optimal upgrade probability of type $n-2$ demand is $\Dstar_{n-2} = \text{argmax}_{v_{n-2}} w^{\textsc{U}}(c_{n-2},\left(c_{n-1}-\frac{\Lambda_{n-1}}{1-\Dstar_{n-1}}\right)^{+},\Lambda_{n-2},0)$. By math induction, we have for $i \in \{n-3, n-4, \ldots, 1\}$, $\Dstar_{i} = \text{argmax}_{v_{i}} w^{\textsc{U}}(c_{i},\left(c_{i+1}-\frac{\Lambda_{i+1}}{1-\Dstar_{i+1}}\right)^{+},\Lambda_{i},0)$. 
\Halmos
\endproof

\proof{Proof of Theorem \ref{thm:ntype}}
Recall that in Theorem \ref{thm:sameopt}, we have defined the upper bound of total revenue $W^{\textsc{U}}$, which directly depends on $c_1$, $c_2$, $\Lambda_1$, and $\Lambda_2$. In addition, it maximizes over all possible $v_1$. Therefore, in this proof, we denote the total revenue of the deterministic programming of two types of resource as $w^{\textsc{U}}(c_1,c_2,\Lambda_1,\Lambda_2;\Dstar_1)$. That is, if $\Lambda_1 = \min\{\Lambda_1, \frac{c_2-\Lambda_2}{\Dstar_1}, \frac{c_1}{1-\Dstar_1} \}$, 
\[
w^{\textsc{U}}(c_1,c_2,\Lambda_1,\Lambda_2;\Dstar_1) = \Lambda_2r_2 + \Lambda_1T \Big( R_1(\Dstar_1)+r_1 \Big),
\]
if $\frac{c_2-\Lambda_2}{\Dstar_1}= \min\{\Lambda_1, \frac{c_2-\Lambda_2}{\Dstar_1}, \frac{c_1}{1-\Dstar_1} \}$, 
\[
w^{\textsc{U}}(c_1,c_2,\Lambda_1,\Lambda_2;\Dstar_1) = \Lambda_2r_2 +\frac{c_2-\Lambda_2}{\Dstar_1} R_1(\Dstar_1) + r_1\min\Big\{\Lambda_1, c_1+c_2-\Lambda_2 \Big\},
\]
if $\frac{c_1}{1-\Dstar_1}= \min\{\Lambda_1, \frac{c_2-\Lambda_2}{\Dstar_1}, \frac{c_1}{1-\Dstar_1} \}$, 
\[
w^{\textsc{U}}(c_1,c_2,\Lambda_1,\Lambda_2;\Dstar_1) = \Lambda_2r_2 + \frac{c_1}{1-\Dstar_1} R_1(\Dstar_1) + r_1\min\Big\{\Lambda_1, c_1+c_2-\Lambda_2 \Big\}.
\]

Let $\mathbf{\Dstar}=[\Dstar_1, \Dstar_2, \ldots, \Dstar_{n-1}, 0]$, where each $\Dstar_i$ is defined in Theorem \ref{thm:upperm}, we have the total revenue of \eqref{eq:optdefm} is 
\[
w^{\textsc{U}} = w^{\textsc{U}}(c_{n-1}^{(1)},c_{n}^{(1)} ,\Lambda_{n-1},\Lambda_n;\Dstar_{n-1})+\sum_{i=1}^{n-2} w^{\textsc{U}}(c_i^{(1)},\left(c_{i+1}^{(1)} - \frac{\Lambda_{i+1}}{1-\Dstar_{i+1}}\right)^{+},\Lambda_i,0;\Dstar_i),
\]
where we input $\left(c_{i+1}^{(1)} - \frac{\Lambda_{i+1}}{1-\Dstar_{i+1}}\right)^{+}$ is because if we take $\Dstar_{i+1}$ as the static upgrade probability, there will be $\frac{\Lambda_{i+1}}{1-\Dstar_{i+1}}$ type $i+1$ demand units not upgraded. 

Define $W^{\pi}_i$ as the expected total revenue generated by Algorithm \ref{alg:ntype} on type $i$ requests. By Theorem \ref{thm:2type}, we have, 
\[
\mathbb{E}\left[w^{\textsc{U}}(c_{n-1}^{(1)},c_{n}^{(1)} ,\Lambda_{n-1},\Lambda_n;\Dstar_{n-1}) - w^{\pi}_{n-1}-w^{\pi}_{n} \right] = O(\log T),
\]
and for any $i \in [n-2]$,
\[
\mathbb{E}\left[w^{\textsc{U}}(c_i^{(1)},\left(c_{i+1}^{(1)} - \frac{\Lambda_{i+1}}{1-\Dstar_{i+1}}\right)^{+},\Lambda_i,0;\Dstar_i) - W^{\pi}_i \right] = O(\log T).
\]

Therefore, we can obtain
\begin{align*}
   & \mathbb{E}\left[ w^{\textsc{U}} -  W^{\pi} \right] \\&= \mathbb{E}\left[w^{\textsc{U}}(c_{n-1}^{(1)},c_{n}^{(1)} ,\Lambda_{n-1},\Lambda_n;\Dstar_{n-1}) - W^{\pi}_{n-1}-W^{\pi}_{n} \right] + \sum_{i=1}^{n-2} \mathbb{E}\left[w^{\textsc{U}}(c_i^{(1)},\left(c_{i+1}^{(1)} - \frac{\Lambda_{i+1}}{1-\Dstar_{i+1}}\right)^{+},\Lambda_i,0;\Dstar_i) - W^{\pi}_i \right] \\&= (n-1) O(\log T) = O(n \log T).
\end{align*}
\Halmos
\endproof

\end{document}